\documentclass[12pt,letterpaper]{amsart}    

\usepackage[parfill]{parskip}    

\usepackage{graphicx}

\usepackage{amssymb} 

\usepackage{amsmath} 

\usepackage{amsthm} 

\usepackage{mathrsfs}

\usepackage{epstopdf}

\usepackage{enumerate} 

\usepackage{fullpage}

\usepackage{hyperref} 

\theoremstyle{plain}
\newtheorem{thm}[equation]{Theorem}

\theoremstyle{definition}

\theoremstyle{remark}

\begin{document}

\title{Lucjan Emil B\"ottcher and his mathematical legacy}

\author[S. Domoradzki]{Stanis\l aw Domoradzki}\footnote{This work was partially supported by the Centre for Innovation and Transfer of Natural Sciences and Engineering Knowledge, University of Rzeszów.}
	
\address{Institute of Mathematics, University of Rzesz\'ow, ul. prof. S. Pigonia 1, 35-959 Rzesz\'ow, Poland}
\email{domoradz@univ.rzeszow.pl}
\author[M. Stawiska]{Ma{\l}gorzata Stawiska}
\address{Mathematical Reviews, 416 Fourth St., Ann Arbor, MI 48103, USA}
\email{stawiska@umich.edu}

\date{\today}                                           

\maketitle


\tableofcontents

\setcounter{tocdepth}{3}

\begin{abstract} This article concerns the life and work of Lucjan (\L ucjan) Emil B\"ottcher (1872-1937), a Polish mathematician. Besides biographical and bibliographical information, it contains a survey of his mathematical achievements in the theory of iteration and holomorphic dynamics. Some documents are presented for the first time. 
\end{abstract}

\section{Introduction}

The name of Lucjan Emil B\"ottcher is familiar to mathematicians interested in functional equations, theory of iterations or dynamics of holomorphic functions. It is associated with B\"ottcher's equation $F(f(z))= [F(z)]^n$, where $f(z) =z^n+a_{n+1}z^{n+1}+..., \quad n \geq 2$, is a known  function which is analytic in some neighborhood of the point $0$ in the complex plane), B\"ottcher's coordinate (the unknown function in B\"ottcher's equation, also called B\"ottcher's function or B\"ottcher's map)  and B\"ottcher's theorem (establishing existence of B\"ottcher's map under certain conditions).  But overall, he remains a relatively obscure figure and his works are very little known.\footnote{He is not mentioned in any of the editions of J. C. Poggendorff's biographical dictionary.} \\

Here we try to present Lucjan Emil B\"ottcher  as a mathematician truly ``without borders". Born in Warsaw (then under the Russian rule) in a Polish family of Evangelical Lutheran denomination (thus a minority among Roman Catholics), he studied mathematics in Warsaw and engineering in Lvov (then in Austro-Hungarian monarchy), took his doctorate in Leipzig and worked as an academic teacher in Lvov (which after 1918 became a part of newly independent Poland).  We write in detail about his education,  people who influenced him, his career and obstacles he faced in its advancement, and his activities in scientific societies. We add  some newly found biographical details to what the first author of this article has already written about B\"ottcher's life and activities (\cite{Do1}, \cite{Do2}).  In particular,  
we compile a complete (to the best of our knowledge) bibliography of B\"ottcher's publications and, on its basis, discuss his mathematical ideas  and results, along with their later development and impact.  We highlight the places where particular concepts appear in his works. For example,  all authors quoting B\"ottcher's theorem (as well as the function and equation named after him) refer to his paper published in    1904 in Russian, while we find out that the result already appeared in a paper of his written in Polish, in 1898.\\

Quite recently, another contribution by B\"ottcher got recognition. Namely, it  was pointed out in \cite{AIR},  \cite{EL} and \cite{Mi2} that  B\"ottcher  gave  examples of rational maps for which the whole sphere is the ``chaotic" set, 20 years before Samuel  Latt\`es independently came up with maps with the same property, which are now known as ``Latt\`es examples".  Although one such  example appeared even earlier in the work of Ernst Schr\"oder, B\"ottcher was  was the first to consider these examples from the dynamical point of view and (according to \cite{Mi2}) he was the first to use the term ``chaotic" in reference to their behavior.  Moreover, as the authors of \cite{AIR} acknowledge, ``he seems to have made a conceptual leap that would not be seen in print in the French study until Fatou's 1906 {\it Comptes Rendus} notice: given a particular function $f$, he viewed the sphere as partitioned into convergence regions by boundary curves" (p. 177). Such a partition arose from the study of the convergence of Newton's method for a quadratic equation undertaken by Arthur Cayley; other  mathematicians, including Schr\"oder, attempted to describe convergence regions in other cases, without much progress.  It was B\"ottcher who  formulated several general properties of boundary  curves of the regions of convergence (the boundaries are now known to be contained in the Julia set) and described them explicitly in the simplest cases (of  monomials and Chebyshev polynomials).  He also stated an upper bound for the number of  ``regions of convergence"  (properly speaking,  of non-repelling cycles) of a rational function in terms of the number of critical points.  The bound was later conjectured by Fatou (who himself proved a  weaker estimate)  and was proved to be sharp, by M.  Shishikura, only in 1980s. However, it should also be noted that B\"ottcher presented very few proofs of his statements and his account is  mostly schematic, sometimes  hypothetic, plainly speculative or even mistaken. This was the main reason for his being underappreciated during his lifetime. His forerunning insights had to wait for full development by other mathematicians who rediscovered them independently. These significant ideas, while few in number,  not only have their place in history, but are very much a part of today's mathematics and deserve to be widely  known. So does their author.  \\

\section{Acknowledgments} We would like to thank the following people and institutions:  Professor  Ludwig Reich (University of Vienna), who inspired the first author's interest in Lucjan Emil B\"ottcher;  also, Dr. Zofia Pawlikowska-Bro\.zek (AGH University of Science and Technology in Krak\'ow), Professor Alexandre Eremenko (Purdue University), Dr. Pawe\l\ Polak (Pontifical Uniwersity of John Paul II  in Krak\'ow), Professor Liliana R. Shakirova (Kazan Federal University), Professor Mykhaylo Zarichnyy (Lviv University).

\section{The life of Lucjan Emil B\"ottcher}

Lucjan Emil B\"ottcher was born on January 7 (21)\footnote{21 is listed in the official documents of the Lvov Polytechnic related to B\"ottcher's retirement; in parentheses there is February 2. In other documents, e.g., B\"ottcher's CV, January 7 is listed. These discrepancies are due to the differences between the Julian and Gregorian calendars. See personal file, L. B\"ottcher, fond 26, op. 5-58, Lvov District Archive.}, 1872, in Warsaw, in a Lutheran family. His father was Piotr, his mother was Anna n\'ee Kraus. In Warsaw,\footnote{There was no Polish state at that time. Poland was partitioned among three occupants: Russia, Austria and Prussia. Warsaw was under Russian occupation, Lvov belonged to the Austro-Hungarian monarchy.} in the years 1881-1885 he attended Herman Benni's \footnote{Herman Benni (1834-1900), an  Evangelical Lutheran  pastor, a graduate of theological studies in Dorpat. In 1880 he opened a private Men's School, which was closed after 5 years at the excuse of lectures being conducted in Polish. This was a time of intensified russification. See p. 239, \cite{Sch}.} four-class real school; later, in 1886-1891, he attended Pankiewicz's \footnote{Jan Pankiewicz (1816-1899), a graduate of the St. Petersburg university, where he obtained the degree of  candidate in philosophy. Since 1841 he taught mathematics in the Real Gymnasium in Warsaw and descriptive geometry in the School of Fine Arts. A school director and inspector, a translator of works in mathematics and chemistry, e.g., A. M. Legendre's {\it Beginnings of Geometry} (1844), an author of entries in mathematics in Orgelbrand's {\it Encyklopedia Powszechna}. The school which B\"ottcher attended was founded in 1876. Pankiewicz managed it until 1894.} six-class real school. Having completed the latter he passed an exam on the material of six classes in a state real school in Warsaw. His intention was to study mathematics at the university  in Warsaw, for which he needed education in classics and the maturity exam (matura). He completed his education in the classical gymnasium in \L om\.za (then in the Polish Kingdom under the Russian rule, which since 1883 was also known as Vistula Country), where he passed his maturity exam in 1893. Then he began studying mathematics at the Imperial University of Warsaw. He was a student in the academic year 1893/94. As he noted in his CV, he had to leave the university because of his participation in a demonstration in honor of Colonel Kili\'nski.\footnote{Jan Kili\'nski, born 1760, a shoemaker by trade, a leader of the burghers of Warsaw during the Ko\'sciuszko insurection in 1794. Imprisoned in St. Petersburg in 1794-1796.} In Warsaw he attended lectures in mathematical analysis (by Nikolay Yakovlevich Sonin), analytic geometry (by Vassily Afanasyevich Anisimov), descriptive astronomy (by Ehrenfeucht), as well as in general chemistry and physics. Later he moved to Lvov and enrolled as a student in the Division of Machine Construction of the Polytechnic School, where he studied in the years 1894/96 - 1896/97. He obtained a so-called half-diploma and passed  the first state exam. As a student in Lvov, he was active in the students' engineering circle. In the c.k. (imperial and royal) Polytechnic School in Lvov a student had to pass the following mathematical subjects in the first two years of study: mathematics- course I (6 hours of lectures in the winter semester and 6 hours of lectures in the summer semester), descriptive geometry (5 hours of lectures and 10 hours of exercises, called repetitory, each semester, respectively), repetitory in elementary mathematics (2 hours each respectively), repetitory in higher mathematics (2 hours each). The programs of these courses are presented in the Appendix.\\

It should be added that for the first state exam in the Division of Machine Construction, the following material in mathematics was compulsory: mathematics, course I and II, and descriptive geometry. It is noteworthy that already in 1895 B\"ottcher published in Lvov lithographed materials for students in differential and integral calculus.\footnote{Information after B\"ottcher's CV. We were not able to find these lectures.}\\

At the beginning of 1897 he interrupted his course of technical studies in Lvov and moved to Leipzig in order to study mathematics.  He  attended there lectures of the following professors: Sophus Lie ({\it Theory of differential invariants}, {\it Theory of differential equations}, {\it Theory of continuous   transformation groups}; seminars {\it Theory of integral invariants} and {\it Differential equations}), Adolph Mayer ({\it Higher analytical mechanics}), Friedrich Engel ({\it Differential equations}, {\it Algebraic equations}, {\it Non-Euclidean geometry}) and Felix Hausdorff ({\it Similarity transformations}). He finished his studies obtaining in 1898 the Doctor of Philosophy degree on the basis of the dissertation ``Beitr\"age zu der Theorie der Iterationsrechnung" (published by Oswald Schmidt, Leipzig) and complying with other procedures. He was extremely industrious; in a very short time he managed to prepare the doctoral dissertation. Formally, his work was supervised by Sophus Lie, an outstanding mathematician, known, among other things, for his study of continuous transformation groups (now called Lie groups).\footnote{Among Lie's students there also were Elie Cartan and Kazimierz \.Zorawski.} B\"ottcher's intention was to treat in his thesis the theory of iteration from the point of view of Lie groups. Even though this ambitious attempt remained largely unsuccessful, he managed to outline several deep ideas. There was some disagreement among the committee members regarding evaluation of his thesis; in the official application to grant B\"ottcher the doctoral degree the signature of Wilhelm Scheibner, who did not want to supply an official report, is struck out. Lie then engaged in correspondence with the university's officials in support of B\"ottcher, explaining importance of B\"ottcher's investigations and their relation with his own research: {\it ``As both the author and Mr. Scheibner indicate the relationship of the submitted work to my concept of one-parameter groups of  transformations, I agree in part with these comments. The relationship, however, lies a little deeper. In 1874, I thought that every finite transformation of a finite continuous group is contained in a one-parameter subgroup. In 1883, I formulated the question whether this fundamental theorem also applies to  infinite continuous groups. However, since this issue  exceeds not only my strength, but also the strength of the current analysis, I restricted myself mainly to showing only for specific examples that this question can be answered in the affirmative.\\ 
Various authors, including Mr. B\"ottcher, considered the same issue for a particularly important group, namely the group of all point transformations. However, I cannot admit that the author has managed to definitively substantiate significant new results. Despite all of this, his considerations, which testify to the diligence  and talent, have their value. (...)
In any case,  I (as well as Mr. Scheibner) agree that this attempt be accepted as a thesis and we also agree regarding the grade II. I am choosing such a good grade  because Mr B\"ottcher himself chose the topic and developed it independently.}\footnote{cf. \cite{Do2}}  Ultimately, B\"ottcher's thesis was evaluated as IIa (admodum laudabilis, the second highest grade) by Lie, who wrote that "the candidate is an intelligent mathematician, possessing good and solid knowledge."  The materials from the University Archive in Leipzig (Universit\"atsarchiv Leipzig, Phil. Fak., Prom., 714, Bl 7) are presented in the Appendix.\\

The Polish version of B\"ottcher's thesis, modified and expanded, was published in three parts in the years 1899-1902 in "Prace Matematyczno-Fizyczne", practically the first Polish journal devoted solely to mathematical and natural sciences, which was privately financed by  a distinguished publisher and historian of mathematics, Samuel Dickstein (1851-1939).\\

After returning from Leipzig to Lvov, B\"ottcher was appointed to the post of an assistant in c.k. (imperial and royal) Polytechnic School in Lvov. He worked there in the period of 1/10/1898 - 30/9/1910, initially in the Chair of Mechanical technology, later (since 1899) in the Chair of mathematics. He had his PhD diploma from Leipzig nostrified in Lvov in 1901. Since that year, he made multiple attempts to obtain license to lecture at the Lvov University. We present the copies of the documents in the Appendix. Let us outline his application for the license to lecture ({\it venia legendi} and habilitation).  He first applied for admission to habilitation on October 2, 1901. The committee appointed for his habilitation proceeding met on on January 13, 1902, and the unanimously negative motion by the committee "not to admit to the further stages of habilitation" was accepted by the Faculty Council on March 14, 1902. In 1911 B\"ottcher obtained {\it veniam legendi} in mathematics in the c.k. Polytechnic School in Lvov. It was possible to obtain habilitation in the Polytechnic School in Lvov since the late eighties of 19th century. This possibility was open also for practicing professionals without PhD. Some of them complemented this procedure later with PhD degrees from other universities and polytechnic schools of the Austro-Hungarian monarchy. We would like to point out  different roles   played by habilitation at universities and at polytechnic schools.\\

Since 1910 B\"ottcher had a position of an {\it adiunkt}; since 1911 he was also a privatdozent of mathematics (additionally, he could lecture in the years 1912 -1918). In the years 1920- 1935, in the period of the Second Republic, he was a docent in the Chair of Mathematics, still employed as an {\it adiunkt}. Before obtaining habilitation at the Polytechnic School he was an assistant to P. Dziwi\'nski. They jointly ran recitations in elementary mathematics and in mathematics course II. B\"ottcher lectured on elements of higher mathematics in the Division of Architecture in 1910- 1927. Between two wars, he lectured in the Division of Civil and Hydraulic Engineering on applied mathematics, theory of vectors, difference equations, notions and methods of elementary mathematics, methods of computation, and, in the Division  of Mechanics, on theoretical mechanics and calculus of variations. In the Appendix we present syllabi of some of his lectures.\\

In 1911 B\"ottcher requested at the Faculty of Philosophy of the Lvov University that his license to lecture at the c.k. Polytechnic School be also recognized at the university. His request was denied. Another time he applied for habilitation in 1918. The committee's opinion was that "the scope and character of the research differs greatly from current scientific mainstream". We present fragments of the committee's decision. Once again, already in free Poland, B\"ottcher requested recognition of his license to lecture on May 1, 1919, taking into account his research work (he was an author of about 20 works in mathematics, see the list and the analysis in further parts of the article). This time, too, the decision was unfavourable.\\

B\"ottcher took part in scientific meetings of mathematicians and philosophers, also in Conventions of Polish Naturalists and Physicians. In many of those conventions there were mathematical and physical sections. He presented results of his research. At  Convention IX (Krak\'ow, 21-24 July, 1900) he presented a talk ``Substitutional functional equations". He considered equations of the type $\Phi\{u(z),u(z_1),z\}=0$, where $\Phi$ is an algebraic rational entire function, $u(z)$ is an unknown function, and $z_1=f(z)$ is some known algebraic function. It should be added that Convention IX met in 22 sections. In the mathematical-physical section (with astronomy) talks were presented by S. Dickstein, J. Puzyna, K. \.Zorawski, and physicists M. Smoluchowski, W\l. Natanson, M. Rudzki. Therefore B\"ottcher was known in the wider scientific community. At Convention X (Lvov, 22-25 July, 1904)\footnote{The convention was planned for the year 1903, but in the same year there were  conventions of physicians in Madrid and Cairo and geological convention in Vienna} he presented a talk ``In the area of theory of functional equations", in which he surveyed some types of functional equations  considered at the previous convention. An eminent mathematician from Krak\'ow, S. Zaremba, also took part in this convention. B\"ottcher also participated in Polish Philosophical Convention I in Lvov in 1923. His talk at this convention, ``On Russell's antinomy" (published in 1927), was cited in W. Tatarkiewicz's ``History of Philosophy", used by many generations of Polish students.\\

B\"ottcher belonged to Polish Mathematical Society in Lvov which was created in 1917 at the initiative of Lvov mathematicians J. Puzyna (1856-1919), Z. Janiszewski (1888-1920), H. Steinhaus (1887- 1972), A. \L omnicki (1881- 1941), P. Dziwi\'nski (1851-1936) and the philosopher T. Cze\.zowski (1889- 1981).\footnote{In 1919, Mathematical Society in Krak\'ow was created, which in 1920 was transformed into Polish Mathematical Society. The operations of the Polish-Bolshevik war interrupted activities of Polish Mathematical Society in Lvov in the second half of 1920. At the proposal of Polish Mathematical Society in Krak\'ow the Society in Lvov dissolved and then reconstituted as the Lvov Branch of Polish Mathematical Society}  It started its activities on December 3, 1918, and its statute was confirmed by the Government's decree number L. XIII a. 30315/452. B\"ottcher is not mentioned among speakers, but he is listed as an active member. Participation in conventions, Society's work, activities in an engineering circle as a student allow us to view him as an extremely engaged and open scientist, devoted not only to mathematics.\\

 B\"ottcher published about 20 papers in significant Polish, Russian and French journals. Their topics and results are presented in a subsequent part of the article. He  also dealt with mechanics: he published a lithographed ``Lecture on general mechanics" and the article ``Some remarks on the principle of inertia". In the (contemporaneous) article \cite{Ku}, B\"ottcher's lecture is characterized as follows: ``[it is]  very carefully developed, including kinematics, statics and dynamics, with examples taking into account the needs of technologists. One can regret that this course was not prepared by the author to be printed as a textbook." As noticed in Chapter 3 of \cite{Pol} , B\"ottcher also made an original attempt to express the foundations of mechanics through some psychological concepts. Additionally, he was concerned  with problems of mathematical education in high school and he wrote high school textbooks. \footnote{He was also involved in proofreading other authors' texts; \cite{Wan} mentions his work on a text by Placyd Dziwi\'nski.} In 1911 he published in Warsaw ``Principles of elementary algebra, adapted to the curriculum in the Polish Kingdom". This textbook was written in the spirit of so-called Meran program, which recommended initiating students into thinking in terms of functions. He wrote: ``Adding to the word `polynomial' the phrase `or an entire function', adding to the words `fractional expression' the phrase `or a rational function' et., and doing so consistently, allows one to familiarize oneself with the notion of function, at least with the simplest one, with which elementary algebra deals." The textbook contains many difficult problems. Let us quote one of them, whose solution allows a student to use knowledge of quadratic functions: ``Two couriers departed from  Skierniewice and from \L owicz  to meet each other. The first one, having left Skierniewice 25 minutes earlier, arrived at \L owicz at the same time at which the other one arrived at Skierniewice. The other passed the first one after 2 hours and 24 minutes of his travel. Knowing that the distance from Skierniewice to \L owicz is 23400 m. compute the time in which each courier traveled his path, and how fast he was going." Note also that the author introduced Lille's graphic method of solving quadratic equations and informs about the ways of using the sliding rule. The text was written with students interested in mathematics in mind. It is vast, with 704 pages. B\"ottcher was also a wholehearted advocate of introducing differential and integral calculus at school. He considered the process of differentiation to be easier than division or taking roots of real numbers, or computing compound interest. The idea of introducing differential and integral calculus occurred  at Convention X of Polish Naturalists and Physicians, after meetings of the section of mathematics and physics and the section of scientific education.  B\"ottcher noticed in talks about education in physics that the ideas of modern physics were contained in school curricula. He also wrote a geometry textbook for high schools, ``Principles of geometry with numerous exercises" (Warsaw 1908). The text contains material in planimetry and has 322 pages. In the reviewer's opinion,\footnote{The review by A. Wilk, \cite{Wi}} the author was going beyond the curriculum, too often relating to the material in elementary physics, but his  drawings were prepared very meticulously, which can facilitate understanding of the material by students. Let us note that B\"ottcher published texts for Polish schools under the Russian rule, despite well-developed  market  for publications in Lvov. Therefore the reviewer wrote: ``It could also be advantageously used in our high schools, despite being somewhat too vast, which however can be amended; it could be particularly recommended to more diligent students for studying and working out numerous exercises."\\

After 1912 B\"ottcher developed interest in spiritualism and metapsychology.\footnote{The interest might have been influenced by Wilhelm Wundt, whose lectures in psychology B\"ottcher attended as a student in Leipzig. The authors thank Professor Cezary Doma\'nski for pointing out this possible connection. One of Wundt's students in 1870's was Julian Ochorowicz, a Polish scientist and inventor, who studied hypnosis, occultism, spiritualism, and telepathy. It should be noted, however, that Wundt was an experimentalist and considered the notion of soul to be irrelevant for explaining psychological phenomena.} He wrote books on these topics, e.g. {\it ``Table-turning"} (1912)\footnote{In the series: Biblioteka Wiedzy Og\'olnej, published by Kultura i Sztuka, Lw\'ow, 1912. There is a copy in the special collection of the Jagiellonian Library in Krak\'ow.}  which is ``a practical outline containing instructions for conducting seances with table-turning".\\

He retired on August 31, 1935.\\

B\"ottcher   married  Maria Wolle (in 1900) and had four children: W\l adys\l aw (born 1901), a Defender of Lvov in 1920 \footnote{cf. Z. Pop\l awski, {\it Dzieje Politechniki Lwowskiej 1844-1945}, Wroc\l aw, 1992}; Marian (1911); and two other, whose  names  are  unknown. Here is how Pawe\l \ B\"ottcher, Lucjan's grandson, remembers his family: {\it ``Both my maternal and paternal grandparents came from Lvov. After the war they were repatriated to Bielsko-Bia\l a. They hoped it was only temporary. They did not want to take root there. They did not talk a lot about what they lost. This is why I know so little today. My paternal grandfather, Lucjan B\"ottcher, taught mathematics at the Lvov Polytechnic, where also my maternal grandfather was  a student..."} Pawe\l 's father was Marian. Mr. Pawe\l \  B\"ottcher, whom the first author contacted, did not have any more information related to his family, as he himself noted in his memoirs from the Eastern Borderlands. \\

Lucjan B\"ottcher died on May 29, 1937, in Lvov.

\section{The bibliography of L. B\"ottcher}

The bibliography was compiled as a result of search in: the journals {\it Jahrbuch Fortschritte der Mathematik}, {\it Prace Matematyczno-Fizyczne}, {\it Wiadomo\'sci Matematyczne}, {\it Muzeum}, the database Zentralblatt MATH, the bibliographic guide by Byelous (items 2-8 and 11 in the list of other publications) and B\"ottcher's personal file in the District Archive in Lviv. The list may be incomplete. The author self-published  some of his works.

\subsection{Mathematical publications}

\begin{enumerate}

 \item  Zasadnicze podstawy teoryi iteracyi [Essential foundations of the theory of iterations], pp. 8, {\it Pami\c etnik Towarzystwa Politechnicznego we Lwowie}, issue  1, s. 126-133. a supplement to {\it Czasopismo Techniczne}, Lw\'ow, 1897. 

\item Beitr\"age zu der Theorie der Iterationsrechnung,  published by Oswald Schmidt, Leipzig, pp.78, 1898 (doctoral dissertation).

\item Przyczynki do teoryi rachunku iteracyjnego [Contributions to the theory of iterational calculus], {\it Wiadomo\'sci Matematyczne}, vol. II(1898), s. 224-229 ( the author's discussion of his doctoral dissertation).

\item Kilka s\l \'ow z dziedziny rachunku iteracyjnego [A few words in  the subject of iterational calculus], an offprint from {\it Czasopismo Techniczne}, vol. 17 (1899), pp.56-57, Lw\'ow.

\item Zasady rachunku iteracyjnego (cz\c e\'s\'c pierwsza i cz\c e\'s\'c druga) [Principles of iterational calculus (part one and two)], {\it Prace Matematyczno - Fizyczne}, vol.  X (1899 - 1900), pp. 65 - 86, 86-101. 

\item R\'ownania funkcyjne podstawnicze [Substitutional functional equations] , {\it Wiadomo\'sci Matematyczne}, vol. IV(1900), s. 233-235

\item O wlasno\'sciach wyznacznik\'ow funkcyjnych [On properties of  functional determinants], {\it Rozprawy Wydzia\l u Matematyczno - Przyrodniczego Akademii Umiej\c etnosci w Krakowie}, v. 38 (general volume) (1901); series II, v. 18(1901), 382-389.

\item Zasady rachunku iteracyjnego (cz\c e\'s\'c III) [Principles of iterational calculus (part III)], {\it Prace Matematyczno - Fizyczne}, v. XII(1901), p. 95-111
 
\item  Zasady rachunku iteracyjnego (cz\c e\'s\'c III, doko\'nczenie) [Principles of iterational calculus (part III, completion)], {\it Prace Matematyczno - Fizyczne}, v. XIII(1902), pp. 353-371 

\item Rozwijanie na szeregi potegowe funkcyi, okre\'slonej r\'ownaniem algebraicznym nieprzywiedlnem f(x,y) = 0 [Expansion into power series of functions defined by an algebraic irreducible equation f(x,y)=0], {\it Wiadomo\'sci Matematyczne}, tom VII(1903), s. 1-21. 

\item Glavn"yshiye zakony skhodimosti iteratsiy i ikh prilozheniya k" analizu [The principal laws of  convergence of iterates and their application to analysis], {\it Bulletin de la Societe Physico-Mathematique de Kasan}, tome XIII(1, 1903), p.1-37, XIV(2, 1904), p. 155-200, XIV(3, 1904), p. 201-234.

\item Iteracye funkcyi liniowej [Iterations of a linear function], {\it Wiadomo\'sci Matematyczne}, vol. VIII(1904), s. 291 - 307.

\item Iteracye funkcyi liniowej (ciag dalszy i doko\'nczenie) [Iterations of a linear function (continuation and completion)], {\it Wiadomo\'sci Matematyczne}, vol. IX (1905), p. 77-86.

\item Notatka o rozwi\c azaniu r\'ownania  funkcyjnego $\Psi f(z) -\Psi(z)=F(z)$  [A note on solution of the functional equation $\Psi f(z) -\Psi(z)=F(z)$], Wiadomo\'sci Matematyczne, tom XIII(1909), s. 17-21. 

\item Nouvelle m\'ethode d'int\'egration d'un syst\`eme de $n$ \'equations fonctionnelles lin\'eaires du premier ordre de la forme $U_i(z)=\sum_{j=1}^{j=n}A_{i,j}(z)U_jF(z)$, {\it Annales scientifiques de l'Ecole Normale Sup\'erieure},  tome 26(1909), p. 519-543.

\item Iteracye funkcyi kwadratowej [Iterations of a quadratic function], {\it Wiadomo\'sci Matematyczne}, v.. XVIII(1914), s. 83 - 132.

\item Zasady rachunku iteracyjnego [Principles of iterational calculus], {\it Wektor} 9(1912), p. 501-513.

\item  Iteratsiya $f_x(z)$ algebraicheskoy funktsii $f(z)$ kak" metatranstsendental'naya funktsiya ukazatelya $x$ [Iteration $f_x(z)$ of an algebraic function as a metatranscendental function of the exponent $x$], {\it Bulletin de la Societe Physico-Mathematique de Kasan}, tome XVIII(3, 1912), p.  106-125.

\item Przyczynek do rachunku iteracyj funkcyi algebraicznej wymiernej ca\l kowitej [A contribution to the calculus of iterations of an algebraic rational entire function], {\it Wiadomo\'sci Matematyczne} XVI(1912), s. 201-206.  

\end{enumerate}

\subsection{Other publications of L. B\"ottcher}

Textbooks in mathematics and mechanics, remarks on mathematics education, mechanics, philosophy, logic, occult/spiritualism\\

\begin{enumerate}

\item Repetytoryum Wy\.zszej Matematyki, Rachunek r\'o\.zniczkowy (55 ss.), Rachunek ca\l kowy (47 ss.),  (litografowany podr\c ecznik), [A review course in higher mathematics, Differential calculus (55 pp.), Integral calculus (48 pp.), (a lithographed text)], Lw\'ow, 1895.

\item O podziale k\c ata na trzy cz\c e\'sci r\'owne [On dividing an angle into three equal parts], {\it Czasopismo Techniczne}, vol. 13(1895), s.106

\item Znu\.zenie w szkole [Weariness at school],  {\it Czasopismo Techniczne}, vol. 17(1899), pp. 4-5.  

\item  Teorya wyznacznik\'ow, cz. 1 [Theory of determinants, part 1] , Lw\'ow, 1899, pp. 118 

\item Z teoryi r\'ownan liczebnych [From the theory of numerical equations], {\it Czasopismo Techniczne}, 18(1900), p. 306-307. 

\item Praktyczne rozwi\c azywanie liczebnych algebraicznych r\'ownan stopni wy\.zszych [Practical solutions of numerical algebraic equations of higher degrees], {\it Czasopismo Techniczne}  19(1901), 7-8, s. 15-17. 

\item Przybli\.zony rachunek rzeczywistych pierwiastk\'ow [Approximate computation of real roots], {\it Czasopismo Techniczne} 19(1901), pp. 91-92, s. 114 

\item Kilka uwag z powodu artyku\l u p. Bronis\l awa Biegeleisena "U podstaw mechaniki", [Some remarks on the occasion of the article by Mr. Bronis\l aw Biegeleisen, "At the foundations of mechanics"], {\it Czasopismo Techniczne} 20(1902), p. 147-148.

\item Obliczanie funkcyj trygonometrycznych kat\'ow bardzo ma\l ych [Computing trigonometric functions of very small angles], {\it Czasopismo Techniczne} 20(1902), p. 255-256.  

\item  
Wyklady mechaniki. Wyd. Kazimierz Bartel, , rps. 
powiel. [Lectures in mechanics. Publ. by K. Bartel \footnote{Kazimierz Bartel, later a professor of descriptive geometry in the Lvov Polytechnic and a prime minister of the government of Poland.}, mimeographed manuscript] Lw\'ow 1904, cm 25, s. 308 (in the collection of the National Library in Warsaw; no title page, table of contents or bibliography)

\item Wyklad mechaniki og\'olnej, (litografowany podr\c ecznik) [A lecture in general mechanics (a lithographed text)], published by Wac\l aw Kuty\l owski-Sok\'o\l\, in quarto, pp. 738  (Lw\'ow 1905) (information after  \cite{Ku})

\item  \'Sp. W\l adys\l aw Folkierski [Late W\l adys\l aw Folkierski],{\it Czasopismo Techniczne} 22(1904), s. 217-221.

\item Kilka uwag o zasadzie bezwladno\'sci [Some remarks on the principle of inertia], {\it Czasopismo Techniczne}, vol. 23(1905), pp. 237-240, 253-255, 269-271.

\item Kilka uwag w sprawie reformy nauczania matematyki w szkolach \'srednich [Some remarks on the reform of mathematics education in high school], Muzeum, vol. 23(1906), pp. 163-168. 

\item Nauka matematyki w szkole \'sredniej [The study of mathematics in high school], {\it Wszech\'swiat},  vol. 26 (1907), pp. 545-549.

\item Zasady geometryi elementarnej do szk\'o\l \ z licznemi \'cwiczeniami (podr\c ecznik do gimnazjum)[Principles of elementary geometry for schools with numerous exercises (a textbook for middle schools)],  publ. M. Arct, Warszawa, 1908.

\item Rektyfikacja elipsy, Rachunkowe i wykre\'slne metody przybli\.zonego wyznaczania obwodu elipsy, [Rectification of an ellipse, Computational and graphical methods of approximate determination of the circumference of an ellipse], {\it Czasopismo Techniczne}, Lw\'ow, t. 26 (1908), s. 185-187, s. 200-203.


\item Zasady algebry elementarnej. Podr\c ecznik i zbi\'or zada\'n dla szk\'o\l , opracowany wed\l ug najnowszych wymaga\'n pedagogicznych, (704 ss.) [Principles of elementary algebra. A textbook and exercise set, compiled according to the latest pedagogical requirements (704 pp.)], Warszawa 1911.

\item Stoliki wirujace [Table-turning], series Biblioteka Wiedzy Og\'olnej [Library of General Knowledge], publ. Kultura i Sztuka, Lw\'ow, 1915, wyd.  II przejrzane i uzupe\l nione (II edition, revised and completed),  Lw\'ow - Przemy\'sl, 1926. 

\item Problemat \.zycia pozagrobowego. Nie\'smiertelno\'s\'c duszy [The problem of afterlife. Immortality of  soul.], publ. Kultura i Sztuka, Lw\'ow, 1915. 

\item O zasadzie sprzeczno\'sci, [On the principle of contradiction], {\it Przegl\c ad Filozoficzny}, vol. 30, issue 4 (1927)), p. 284. 

\item O antynomji Russela [On Russell's antinomy], {\it Przegl\c ad Filozoficzny}, 30, issue 4 (1927), p. 291-292.

\end{enumerate}

\section{Mathematics of Lucjan Emil B\"ottcher}

\subsection{Glossary}

In order to describe in some detail B\"ottcher's contributions to the area of iterations of rational functions on the Riemann sphere (which is the classical core of complex dynamics in one variable), let us start with a glossary of basic notions and properties related to his  results (following \cite{Ad}, \cite{Mi1}, \cite{Be}, \cite{HY}, \cite{KH}). \\

Let $M$ be a set and let $f: M \mapsto M$ be a mapping. We define the iterates of $f$ as $f^0=Id$ and $f^{n+1}=f\circ f^n$, $n=0,1,...$. The {\bf orbit} of a point $p$ is the set $\{f^n(p): n \geq 0\}$.  A point $p \in M$ is {\bf periodic} for $f$ if there is an $m \geq 1$ such that $f^m(p)=p$. In particular, fixed points are periodic. An orbit of a periodic point is also called a {\bf periodic cycle}. A point is called preperiodic if there are some $m,n >0$ such that $f^m(p) =f^{m+n}(p)$. The {\bf backward orbit} of a point $p$ is the set $\bigcup_{n \geq 0}f^{-n}(\{p\})$ (the union of all preimages of $p$ under all iterates of $f$).\\
Periodic cycles are examples of sets which are (forward) invariant under $f$. More precisely, a set $E \subset M$ is {\bf forward invariant} under $f$ if $f(E)=E$, and it is {\bf backward invariant} under $f$ if $f^{-1}(E)=E$.\\

When $M$ is a topological space, one can study behavior of iterates of a continuous map $f$ near its fixed points. A fixed point $p$ of $f$ is called {\bf topologically attracting} if it has a neighborhood $U$ such that all iterates $f^n$ are defined in $U$ and the sequence $f^n|_U$ converges uniformly to the constant map with value $p$ in $U$. The {\bf basin of attraction} of $p$ is the set of all points $z$ such that the sequence $z, f(z), f^2(z),...$ converges to $p$. A fixed point $p$ is {\bf topologically repelling} if it has a neighborhood $U$ such that for every $p' \neq p$ in $U$ there exists an $n \geq 1$ such that $f^n(p')$ lies outside of $U$.\\
A map is {\bf topologically transitive} if for every pair of non-empty open sets $U$ and $V$ in $M$, there is a non-negative integer $n$ such that $f^n(U)\cap V \neq \emptyset$.\\
Two maps $f$ and $g$ on $M$ are {\bf topologically conjugate} if there exists a homeomorphism $h:M \mapsto M$ such that $h^{-1} \circ f \circ h = g$.\\

Let $(M,d)$ and $(N,d')$ be metric spaces and assume that $N$ is compact. A family $\mathcal{F} \subset \mathcal{C}(M,N)$ is  {\bf normal} if every infinite sequence in $\mathcal{F}$ contains a subsequence  which converges uniformly on every compact subset of $M$. \\

Consider now the case when $M=N=\hat{\mathbb{C}}$, the Riemann sphere (with e.g. the spherical metric), and
$f: \hat{\mathbb{C}} \mapsto \hat{\mathbb{C}}$ is a rational map of degree $d  \geq 2$, i. e., $f=P/Q$ for coprime complex polynomials $P,Q$ with $\max(\mbox{deg }P, \mbox{deg }Q)=d \geq 2$. The properties of a fixed point $p$ of $f$ are determined by its so-called multiplier, i.e., the value $f'(p)$ when $p \in \mathbb{C}$; the multiplier at $\infty$ is understood as the derivative of $1/f(1/z)$ at $z=0$. We say that $p$ is superattracting if $f'(p)=0$, attracting if $0<|f'(p)|<1$; repelling  if $|f'(p)|>1$; rationally neutral if $f'(p)$ is a root of unity (called also parabolic if $f'(p)=1$); irrationally neutral if $|f'(p)|=1$ but $f'(p)$ is not a root of unity. It can be proved that topologically attracting (resp. repelling) fixed points are exactly those that are attracting or superattracting (resp. repelling) in the sense of the multipliers. By the chain rule, the same classification using multipliers can be applied to periodic orbits. \\

The notion of a normal family plays a major role in the study of dynamics of  rational maps.  There are several criteria for the family of holomorphic maps on a subset of the sphere to be normal. For example, a theorem by P. Montel (who introduced the notion) states that a uniformly bounded family of holomorphic functions on a domain in $\mathbb{C}$ is normal.  Let us also mention F. Marty's criterion: A family  $\mathcal{F}$ of holomorphic maps on a region $U \subset \hat{\mathbb{C}}$ is normal if and only if for every compact $K \subset U$ there exists a constant $C_K$ such that $\frac{|f'(z)|}{1+|f(z)|} \leq C_K$ for all $z \in K$ and all $f \in \mathcal{F}$.The quantity on the left-hand side of the inequality is called the spherical derivative of $f$.\\

Using the notion of a  normal family, one can divide the sphere into two subsets. The {\bf Fatou set} is the maximal open subset of $\hat{\mathbb{C}}$ on which the family of all iterates of $f$ is normal and the {\bf Julia set} is its complement. The Julia set is nonempty when $d \geq 2$, but the Fatou set may be empty (as demonstrated by Latt\`es examples). An example of a map with non-empty Fatou set is $f(z)=z^2$: it is straightforward to show that its Julia set is the unit circle. It is not always easy to determine Julia sets, but there are ways to approximate them in the Riemann sphere and to produce computer pictures of them. One such way is to look at the backward orbit of a "typical" point. Namely, if the backward orbit of a point $a$ is infinite (which happens for all but at most two points in the sphere), then its closure contains the Julia set of $f$.  Moreover, if $a$ is in the Julia set of $f$, then the closure of its backward orbit equals the Julia set.\\

Attracting and superattracting cycles for $f$ belong to the Fatou set, while the repelling cycles belong to the Julia set.  In fact, repelling cycles are dense in the Julia set. In general, periodic cycles in the Fatou set (if any) can be identified by considering limit functions of the sequences of iterates of $f$. By Sullivan's theorem, every component of the Fatou set is preperiodic. Assume that the component $F_0$ is  forward invariant. First consider the situation in which every limit of a subsequence of $f^n$ is constant. Then all these constant functions are equal to $\zeta$ such that $f(\zeta)=\zeta$, and $f^n \to \zeta$ locally uniformly on $F_0$.  Hence $\zeta$ is an attracting or superattracting fixed point. If there are non-constant limit functions, then the identity is among them, and the map $f$ is analytically conjugate to an irrational rotation of the unit disk or an annulus (depending on whether $F_0$ is simply or doubly connected). One can also say something when there is a point $ \zeta \in \partial F_0$ (hence in the Julia set) such that $f^n \to \zeta$ locally uniformly as $n \to \infty$. Then necessarily $f'(\zeta)=1$.\\

For a rational function $f$ of degree $d \geq 2$ with a neutral fixed point $z_0$, $f$ is locally conjugate to its linear part if and only if $z_0$ is in the Fatou set. For $f$ with an irrationally neutral fixed point, $f'(z_0)=e^{2\pi i \xi}$, $\xi$ an irrational real number, there are sufficient and necessary conditions of linearizability formulated in terms of approximation of $\xi$ by rationals (given by H. Cremer, C. L. Siegel, A.D. Bryuno, and J.-C. Yoccoz), but much still remains to be investigated.\\ 

 By using perturbation methods it can be proved that a rational map of degree $d$ can have at most $2d-2$ non-repelling cycles. This upper bound equals the number of critical points of $f$, which are the zeros and poles of $f'$. The estimate is sharp, but this was proved only in 1980's by M. Shishikura (\cite{Shi}).\\

 One can construct maps having only repelling  cycles. The {\bf Latt\`es example} gives such a map satisfying $\mathcal{P}(2z)=f(\mathcal{P}(z))$, where $\mathcal{P}$ is the Weierstrass elliptic function associated with a certain lattice in $\mathbb{C}^2$. The construction proceeds as follows (for the details, see \cite{Be}, Section 4.3). A lattice $\Lambda$ is defined as $\Lambda =\{m\lambda +n\mu: m,n \in \mathbb{C}\}$, where $\lambda$ and $\mu$ are fixed complex numbers which are not real multiples of each other. One says that a non-constant function $f$ is  elliptic  for $\Lambda$ if it is a meromorphic function on $\mathbb{C}$ and each $\omega \in \Lambda$ is a period of $f$. E.g., the Weierstrass elliptic function $\mathcal{P}$ is defined as $\mathcal{P}(z) = \frac{1}{z^2}+\sum^*[\frac{1}{(z+\omega^2)}-\frac{1}{\omega^2}]$, where $\sum^*$ is taken over all non-zero elements in $\Lambda$. For $\mathcal{P}$ the addition formula $\mathcal{P}(2z)=R(\mathcal{P}(z))$ holds with the rational function $R(z)=\frac{z^4+g_2z^2/2+2g_3z+(g_2/4)^2}{4z^3-g_2z-g_3}$, where $g_2,g_3$ are uniquely defined in terms of $\Lambda$. This function $R$ has the whole sphere $\hat{\mathbb{C}}$ as its Julia set.  One can show that the family $R^n$ is not equicontinuous on any open subset of $\hat{\mathbb{C}}$: let $D$ be any disk in $\mathbb{C}$ and let $U =\mathcal{P}^{-1}(D)$. For $n$ large enough, the dilation $2^nU$ contains a period parallelogram $\Omega=\{z+s\lambda+t\mu: 0 \leq s,t \leq 1\}$ of $\Lambda$. By the addition formula and periodicity, for these $n$ we have $R^n(D)=R^n(\mathcal{P}(U))=\mathcal{P}(2^nU) = \hat{\mathbb{C}}$, which implies that the Fatou set for $R$ is empty. For a comprehensive treatment of Latt\`es maps, see \cite{Mi2}.\\



Rational maps for which the Julia set is the whole sphere can also be characterized by the following property: $J(R)=\hat{\mathbb{C}}$ if and only if there  is some point $z$ whose forward orbit $\{R^n(z): n \geq 1\}$ is dense in the sphere.\\

B\"ottcher was interested in studying iterations of arbitrary (not necessarily integer) exponent and saw the theory of Lie groups as a suitable framework for this. 
A {\bf Lie group} $G$ is a smooth manifold and a group, with the operations of product $\mu: G \times G \mapsto G$ and inverse $i: G \mapsto G$ such that $\mu$ and $i$ are smooth. An example is $\mathbb{R}^n$ considered as an additive group (with an atlas given by the identity map). \\

 A {\bf one-parameter subgroup} of $G$ is a homomorphism of Lie groups $\theta: \mathbb{R} \mapsto G$ (considering $\mathbb{R}$ as a Lie group), i.e, a homomorphism of groups which is a smooth map. There is one-to-one correspondence between one-parameter subgroups of $G$ and vectors in $G_e$, the tangent space to $G$ at $e$, where $e$ is the unit in $G$. It is given by $\theta \mapsto \theta'(0,1)$, where $(0,1)$ is the unit tangent vector to $\mathbb{R}$ at the origin. This correspondence can be used to define the exponential map $\exp: G_e \mapsto G$ as follows: let $v \in G_e$ and let $\theta_v$ be the corresponding one-parameter subgroup of $G$. Then $\exp(tv):=\theta_v(t)$. The exponential map generates all one-parameter subgroups of $G$ as follows: if $\theta$ is a one-parameter subgroup, then $\frac{\theta(t+h)-\theta(t)}{h}=\theta(t)\frac{\theta(t)-e}{h}$. Hence $\frac{d\theta}{dt}(t)=\frac{d\theta}{dt}(0)\theta(t)$, so $\theta(t) = \exp(t\frac{d\theta}{dt}(0))$.\\


The problem of embedding (a discrete semigroup of) an iterated function into one-parameter group of transformations, which was of interest to B\"ottcher, is in general complicated. Some maps admit such embedding while for others there are obstructions.  See the survey of \cite{BRS} for discussion and references. One of the cases  
for which the question  has been answered in full is that of analytic maps with a rational fixed point. We refer to  \cite{IY}, Chapter IV.  Here we only can say, in very rough terms, that a functional invariant called the Ecalle-Voronin modulus is constructed. 
The analytic classification theorem for parabolic germs states that every parabolic germ is associated with a unique modulus, which is the same for all analytically conjugate germs, and conversely, two formally conjugate parabolic germs with the same modulus are analytically conjugate. 
A parabolic germ can be represented as a flow map of a holomorphic vector field if and only if its Ecalle-Voronin modulus is trivial. The modulus is also used to answer other questions about parabolic germs, e.g., the extraction of iterational roots.\\

A related problem is that of analytic iteration: for a given analytic function $f(z)$  find a  function $\Phi(w,z)$ analytic in $z$ and continuous in $w$ such that $\Phi(w+u,z)=\Phi(w,\Phi(u,z))$ and $\Phi(1,z)=f(z)$.  See \cite{Al} for some historical discussion and references, and \cite{Co}, \cite{EGRS}, \cite{KRS} for the relation between B\"ottcher equation and fractional iteration.\\

\subsection{The works of Lucjan Emil B\"ottcher}

Below we will discuss the content of B\"ottcher's research  in mathematics. We included all publications we could access. \\

\begin{enumerate}

\item ``Contributions to the theory of iterational calculus" (Polish) In this self-report B\"ottcher announces  forthcoming publication of his doctoral thesis in German and in Polish. He mentions its main topics such as a general notion of iteration, its relation with groups of transformations and convergence of iterations.\\

\item B\"ottcher's doctoral thesis      

  In the foreword to his dissertation ``Beitr\"age zu der Theorie der iterationsrechnung..." , written under the supervision of Sophus Lie and published in Leipzig in 1898, B\"ottcher states that his intention is not to give a complete presentation of the theory of iteration, but only of its most important problems. He sees the significance of theory of iteration in its relation with functional equations. He sets three goals for himself: 1) the study of iterated functions with exponents of iteration that are positive integers; 2)  generalizing iteration to admit an arbitrary exponent: negative, rational, real or complex; 3) the study of borderline between 1 and 2 above. The last point is unclear; perhaps after succeeding in theoretically generalizing iteration to admit an arbitrary exponent, he intended to go on to study iteration of functions with these arbitrary exponents. This brings him to seek an interpretation of iteration in terms of groups of transformations,  the theory of which was developed by  Lie. \\
Chapter I of the dissertation is devoted to expressing iterations by means of one-parameter continuous groups of transformations.  The laws of addition, multiplication and association are formulated to illustrate that the collection of iterated transformations with an arbitrary basic transformation forms a group. The analogy  
between arithmetic operations and iterations is highlighted. Iterations with arbitrary exponents are introduced as solutions of certain functional equations or limits of certain expressions. The need to consider complex numbers arises when solving the equation $xa^m=1$ in order to define iterations with negative exponents, or solving $KW(x_1,...,x_n)=WL(x_1,...,x_n)$ to get iterations with fractional exponents. It is observed that the fractional power is not a plain  number, but a ``manifold". Fundamental theorems are formulated on p. 20-21. They are identified as special cases of some results by Lie.  A relation is given between iterations and solutions of some differential equations. \\

Later on, the general considerations give way to the specific study of rational maps of the Riemann sphere. The question of convergence of iterates appears  on p. 33, where the notion of an ``iterational chain" is introduced.   In modern terminology, it corresponds to the notion of an orbit of a point, where the successive images are additionally joined by arcs. It should be noted that the orbits are studied on their own, and various possibilities for limit functions of orbits are investigated, even though at that time the  interest was mainly in the periodic behavior of whole functions rather than single orbits, i.e., in functions satisfying Babbage equation $f^n(z)\equiv z$. B\"ottcher notices that among the rational functions on the Riemann sphere the only functions that are ``periodic" are the M\"obius transforms $z \mapsto (az+b)/(cz+d)$, as the degree of the rational map rises under iteration. The ``points of convergence" are identified as fixed points of the iterated transformation, and it is observed that the magnitude  of $f'(z)$ plays a role in the classification of fixed points. The case $|f'(z)|<1$ gives rise to ``neighborhoods of convergence points", or, as they are known nowadays, basins of attraction. In the case $|f'(z)|=1$ the distinction is made between $\mbox{arg }f'(z)$ being rational or not, but more details are provided only when $f'(z)=1$. It is stated that in the case $f'(z)=1$ the fixed point $z$ lies on the boundary of the ``region of convergence". A description of petals (attracting and repelling) is given, following the work of L. Leau.\\

More interestingly, the boundaries of the regions of convergence (now known to be contained in the Julia set) come into focus. The following theorem is stated on p. 62:  The boundary of the subregion associated with the point $x$ which belongs to the region of convergence, is a closed curve (fully determined in the sphere), which is invariant under the transformation $x_1+iy_1=f(x+iy)$. This theorem is not true in general, as the Julia set of a rational map (which  equals the boundary of the Fatou set, if the latter is not empty) does not have to be a closed curve ( e.g., for $z \mapsto z^2+i$) or even a connected set (e.g., $z \mapsto z^2 + 4$), but it is true that it is totally invariant under the map $f$. Additionally, it is stated that ``the whole iteration chain whose starting point is on the boundary of the subregions of convergence is fully registered by this boundary", which may be interpreted as the statement that the  map $f$ on its Julia set is topologically transitive. In the process of determining the basin of attraction by iterating backward a neighborhood of an attracting point (which mentions preimages of attracting periodic points, called by B\"ottcher ``\.Zorawski points", after Kazimierz \.Zorawski, who studied them in his work on iteration), the way to approximate the ``boundary curve" using backward orbits of points is outlined.\\

Even though the notion of a normal family was not available at that time, some of B\"ottcher's insights can be formulated in terms of normal families of holomorphic maps. First of all, he discussed the regions of convergence (taking into account only basins of attraction and parabolic petals) and of nonconvergence for iterates (he gave examples in which such a region is the whole sphere; nowadays we know that the Julia set which is not the whole sphere must be nowhere dense). Also, one of the few theorems for which he sketches a proof (in section 37) reads as follows: ``A fundamental property  of the invariant curves of the basic function $f(z)$ which simultaneously form the boundary of arbitrary subregions is that $\lim_{n \to \infty} f_n'(z) \neq 0$."  In view of Marty's criterion, it can also be interpreted as the property that on the Julia set the iterates of $f$ do not form a normal family. 

 Section 38 starts with the observation that all  points  of ``negative convergence" (i.e., repelling points) correspond to invariant curves of $f$. On p. 63 there is an example of a rational function without ``regions of convergence", i.e., a function whose Fatou set is empty. It is $f(z)= (Az^4+B(2z^2-1))/(B+A(2z^2-z^4))=\rm{Cn}(2)\rm{Cn}_{-1}(z)$, where $\rm{Cn}$ is the Jacobi elliptic cosine function. Later, Samuel Latt\`es would independently construct  examples of maps whose Julia set is the whole sphere, also starting with elliptic functions. \\

In Part III, B\"ottcher studies a problem of defining an iterate with a general base and a general exponent as a solution of a functional equation.\\

\item ``Principles of iterational calculus" (Polish) There is substantial overlap between this paper and B\"ottcher's doctoral thesis, published almost simultaneously; perhaps this is what he referred to when announcing publication of his PhD thesis in Polish. However, some new content appears in this work. \\
Again, the main idea is to treat iterational calculus as a special chapter of the theory of one-parameter continuous groups of transformations. B\"ottcher revisits the group laws and their use in defining of iteration with an arbitary exponent. He points out to the role of conjugacy in these definitions and the ways to obtain conjugacy as a solution of a functional equation, attributing the point of view to Babbage. A new detail is the construction of an infinitesimal generator of one-parameter group, as well as of the "iterative logarithm" of a transformation $f$, following  Korkine and Farkas.\\
As in the thesis, a good deal of attention is devoted to convergence of iterations of rational functions. Some new light is thrown on the subject. It is made explicit that the investigation of division of the complex plane or sphere into the regions of convergence and divergence of iteration of a rational function is inspired by the work of Arthur  Cayley, who succeeded in completely describing such division coming from the study of Newton-Fourier method for a quadratic equation over complex numbers. ``Fundamental theorems" concern limits of orbits under iteration and ``points of convergence" are classified.  A Latt\`es-type example appears as an illustration of ``transcendental convergence", different from the example presented in the thesis. There are more details about the points which belong to the ``unbounded limit group". In later terminology, these are properties of the Julia set. First of all, this ``group" consists of finite points spread infinitely densely on a certain curve. This, like many B\"ottcher's statements, is not completely clear, and can be interpreted in two ways: the statement says that repelling periodic points are dense in the Julia set, or perhaps it says that the map is topologically transitive on its Julia set. (Both properties hold for rational maps of the sphere.) The "boundary curves of the convergence regions" for some maps are correctly identified: for $f(z)=z^r, r=2,3,...$ the curve is the unit circle, while for $f(z)=2z^2-1$ and for $f(z)=3z-4z^3$ (as well as for higher degree Chebyshev polynomials) the curve is the segment $[-1,1]$ on the real axis.  So B\"ottcher gives first examples of Julia sets. He also presents a new example of an everywhere chaotic rational map 
\[
f(z)=\frac{4z(1-z)(1-k^2z)}{(1-k^2z^2)^2}=\rm{sn^2}2\int_0^{\sqrt{z}}\frac{du}{\sqrt{(1-u^2)(1-k^2u^2)}},
\]
where $\rm{sn}$ denotes the elliptic sine function and the path of integration is arbitrary. An argument is sketched  that this function displays ``transcendental convergence" (i.e., chaotic behavior) in the whole complex plane by noticing that any  limit point of the iterates $f^n$ (taken pointwise) would equal $\rm{sn}^2(\infty)$, which is not a well defined quantity.\\
 \\

Other novel developments are: 1) a conjecture about behavior of an analytic function around an irrationally neutral fixed point, and  2) the upper bound for a number of ``regions of convergence" by the number of critical points of the map $f$ (stated without proof). B\"ottcher  writes in the remark after Section 51: ``It seems that a point $x$ giving $f(x)=x, \ f^{(1)}(x)=e^{\mu \sqrt{-1}}$ with $\mu$ irrational with respect to $\pi$ cannot be a limit point of positive, negative, or even singular convergence. [...]  Whereas it turns out  that the number of all  regions of positive and
singular convergence, regular as well as rhythmic, cannot exceed the number
of the zeros and poles of the derivative of the iterated function, at least
in the case   when the latter is algebraic rational. Therefore the number
of regions of dependence is definitely finite."

 In  this paper B\"ottcher also introduces the equation which now bears his name:

\[
Ff(z)=(F(z))^m,
\]

where $f$ is a given function which is analytic in a neighborhood of its fixed point $x$. 

He proposes a method of solving this equation. When $$f(z)=x+\frac{(z-x)^m}{m!}f^{(m)}(x)+\frac{(z-x)^{m+1}}{(m+1)!}f^{(m+1)}(x)+...,$$ i.e., when the fixed point $x$ is superattracting, he considers what he calls a ``basic algorithm" 
\[
R(z):=\lim_{n \to \infty}\sqrt[m^n]{f_n(z)-z}
\]

and proceeds as follows: Let $\Xi(z)$ be a function such that $\Xi f(z)=\Xi(z)$. Then the general solution to the functional equation is given by
\[
Q^{\Xi(z)(\log R(z))},
\]
where $Q$ is an arbitrary constant. No proof of convergence is given.\\

\item ``Main laws of convergence of iterations and their analytic applications" (Russian)\\

This is the most known and quoted of B\"ottcher's papers. Again, it overlaps  with  the dissertation as well as with ``Zasady". The first part starts with discussion of iteratively periodic functions and the observation that the only rational iteratively periodic functions are of the form $f(z)=(Az+B)/(Cz+D)$. Then limits of iterates of an analytic   functions are considered, the main case being $f^n(z) \to x$ with $f(x)=x$, but it is conjectured (though not investigated) that in some cases the limit of $f^n(z)$ may depend on the point $z$ (i.e., be a non-constant function). The same Latt\`es-type function as in ``Zasady" (involving elliptic sine) is given as an example of ``chaotic convergence". It is stated that for this function any complex number can be obtained as a limit of iterates of $f$ (i.e., in modern terminology, $f$ has orbits that are dense in the sphere). \\
B\"ottcher comes back to the question he had asked in his thesis: that about existence of one-parameter continuous group of transformations containing a given function $f$ and all its iterates.  He admits a possibility that the solution may exists only in some regions of the plane.

In subsequent parts, he analyzes in detail the dynamic behavior of $f(z)=z^2$: the convergence of $f^n$ to $0$ for $|z|<1$ and to $\infty$ for $|z|>1$. He gives a long list of references on ``iterative periodicity" and iteration in general, including Newton's Opuscula. \\

He studies the case of $f$ with an attracting fixed point $x$ and constructs the basin of attraction starting with a disk $C$ such that  $f^n(z) \to x$ in $C$ and taking $\bigcup_nf^{-n}(C)$. \\

Next he considers auxiliary functions $B(z),D(z)$ and $L(z)$ which are used to solve various functional equations. He studies properties of these functions (domain of existence, zeros, poles etc.)\\

Iterative properties of  $f(z) = x+\frac{f^{(m)}(x)}{m!}(z-x)^m+...$ are investigated in detail in Chapter II. B\"ottcher's equation is solved using the solutions to Gr\'evy's equation. First, a function $F$ is defined by the relation $1+F(z)=\frac{(f(z)-x)\bigl(\frac{f^{(m)}(x)}{m!}\bigr)^{-1}}{(z-x)^m}$ and then the infinite product $(1+F(z))(1+F(z_1))...$ is used to find  solutions of the equations  $Yf(z)=\frac{f^{(m)}(x)}{m!}(z-x)^{m-1}Y(z), \quad Y'(x)=1$ as well as $\Psi f(z)=\frac{f'(x)}{m}\Psi(z)$ and $\beta f(z)=m\beta(z)$. Then the function $B(z)$ is defined as $\int_z^x \frac{du}{\Psi (u)} + \frac{1}{m-1}\int_x^z\frac{du}{\Psi (u)}$. It is observed that $B$ has logarithmic singularities at $x$ and all its preimages under $f$. The function $N = e^B$ solves the equation $Nf(z) =[N(z)]^m$.\\

The remaining part of the chapter presents the study of a function with a parabolic fixed point, mainly following Leau.

\item ``On properties of some functional determinants" (Polish) B\"ottcher studies certain functional determinants which are analogous to Wro\'nskians but contain difference operators instead of differential ones.\\

\item ``Rozwijanie na szeregi potegowe funkcji okre\'slonej rownaniem algebraicznem nieprzywiedlnem  $f(x,y)=0$"
This is an expository paper. B\"ottcher discusses expansion of $y$ in (fractional) powers of $x-\bar{x}$, mainly in the case when $\bar{x}$ satisfies the system $f(x,y)=0, \quad f_{0,1}^{(1)}(x,y)=0$. He outlines the procedure based on construction of a certain polygon in the plane, called Newton-Puiseux polygon. He gives references to current texts but not to original works of Newton or Puiseux, and does not mention these names.\\

\item ``Iteration of a linear function" (Polish) B\"ottcher studies iteration of functions of the form $f(z)=(Az+B)/(Cz+D)$. He determines that such a  function is  iteratively periodic (i.e., it satisfies $f^k(z)\equiv z$ for some $k$) if and only if it is an elliptic function, i.e., such that the ratio $(Cr_2+D)/(Cr_1+D)$ is a number of unit modulus, but not $1$. Here  $r_1,r_2$ are fixed points  of $f$.  The period of iteration is rational resp. irrational if $\frac{\arg (Cr_2+D)-\arg(C_1+D)}{2\pi}$ is rational resp. irrational. An example of such a function is $f(z)=\frac{a+bi)z+(c+di)}{(-c+di)z+(a-bi)}$.  Seeking a general formula for iteration of linear functions, he gives recurrence formulas for $A_n,B_n,C_n,D_n$ in  $f_n(z)=(A_nz+B_n)/(C_nz+D_n)$. He also expresses the iterates $f_n$ in terms of symmetric functions of two (distinct) fixed points of $f$ and treats expansion of iterates into continued fractions.\\

\item ``A note on solution of the functional equation $\Psi f(z) -\Psi(z)=F(z)$" (Polish)

In this short note, B\"ottcher  finds a particular solution to the equation mentioned in the title under the assumptions that $f(z)=x+\frac{f'(x)}{1!}(z-x)+..., \quad 0<|f'(x)|<1$ and $F(z)=A_0+\sum_{n=1}^{\infty}[A_{+n}(z-x)^{m/n}+A_{-n}(z-x)^{-m/n}]$ converge in some annulus $r_1<|z-x|<R_1$. Using "Koenigs's algorithm" $B(z)=\lim_{n \to \infty}\frac{f_n(z)-x}{[f'(x)]^n}$ , which solves $Bf(z)=f'(x)B(z)$, he changes the coordinate and expands $F$ in fractional powers of $B$: $F(z) = B_0+\sum_{n=1}^{\infty}[B_nB^{m/n}(z)+B_{-n}B^{-m/n}(z)]$, convergent  in  some $r_2<|B(z)|<R_2$. The function $\psi_0(z)=\frac{B_0 \log B(z)}{\log f'(z)}+\sum_{n=1}^{\infty}\frac{B_n B(z)^{m/n}}{f'(x)^{m/n}-1}+\frac{B_{-n}(B(z))^{-m/n}}{f'(x)^{-m/n}-1}$ is his solution to the title equation.\\

\item ``Nouvelle m\'ethode d'int\'egration d'un syst\`eme de $n$ \'equations fonctionelles lin\'eaires du premier ordre..." (French) 
In this paper, B\"ottcher gives a general solution to the system of $n$ linear functional equations of first order, of the form 
\[
U_i(z) =\sum_{j=1}^nA_{ij}(z)U_j(f(z)), \quad i=1,...,n
\]
using a  solution of the equation $Bf(z)=hB(z), \quad |h|<1$.  He proves the result by constructing a formal Laurent series and investigating its convergence. He considers in detail some special cases of $f$, e.g. $f(z)=az$ or $f(z)=az^b$. Another result of the paper is a ``fundamental law" which concerns the
method of finding a complete and general solution of the mentioned system. The final
part of the paper is devoted to discussion of some properties of the obtained
solution, e.g. its singular points, zeros and points of ramification.\\

\item ``Iteration of a quadratic function" (Polish) B\"ottcher considers iteration of $f(z)=Az^2+2Bz+C$, not necessarily with a positive exponent. He discusses iterations of $f(z)=z^2$ and of $f(z)=2z^2-1$ with an arbitrary exponent by considering them on infinitely sheeted Riemann surfaces. He determines an orbit of an arbitrary point under these iterations. He also introduces the use of ``canonical forms" to aid  iteration of quadratic functions. For  example, he shows that the group of iterations of $f(z)=Az^2+2Bz+C$ can be transformed onto the group of iterations of $\varphi(z)=z^2+T, \quad T=AC-B^2+B$ by $\omega_1\varphi=f\omega_1$ or $\omega_2f=\varphi\omega_2$, where $\omega_1(z)=(z-B)/A, \quad \omega_2(z)=Az+B$.\\

\item ``A contribution to the calculus of iterations of an algebraic rational entire function" (Polish) B\"ottcher considers the following problem: given that $n^{\nu} = p, \quad \nu < 1, \quad p<n$ find whether the entire function $F(z)=A_0z^n+...+A_{n-1}z+A_n$ (i.e., a polynomial) of degree $n$ has a $\nu$-th iterate $F_\nu(z)=a_0z^p+...+a_{p-1}z+a_p$ which is also an entire function, of degree $p$. By expanding $[F(z)]^{p/n}$ in powers of $z^{-1}$ in a neighborhood of infinity, he obtains conditions for solutions (involving certain functions commuting with each other). He also solves an analogous problem for $p > n$.\\

\item ``The notion of iteration..." (Polish) B\"ottcher revisits the problem of defining an iterate with an arbitrary ``base" $f$ and exponent $p$. For this he uses a function $F$ (which he calls an ``iterational functial"[sic]) such that $f(z)=\lim_{n \to \infty}(1+\frac{1}{n}F(z))^{\circ n}$; he takes $f^{\circ p}= \lim_{n \to \infty}(1+\frac{p}{n}F(f))^{\circ n}(z)$ (the superscript $\circ k$ denotes the iteration with exponent $k$). As in his PhD thesis, he claims that such a function $F$ can be constructed in a basin of a point $x$ in the complex plane such that $f(x)=x$ and $0<|f'(x)|<1$ by method of A. Korkine: $F(z)=q\lim_{n \to \infty}\bigl [ \frac{f^{\circ n}(z)-x}{(f^{\circ n})'(x)}\bigr ], \quad q=\log f'(x)$. He proposes construction of such a function in the case $f'(x)=0$ or $1$, but gives no proofs.\\

\item ``Iteration $f_x(z)$ of an algebraic function as a metatranscendental function of the exponent $x$" (Russian)\\

  B\"ottcher studies conditions under which the iteration $f_x(z)$ is a differentiable function of the exponent $x$ (a problem already outlined in his PhD thesis). This brings him to consideration of differential equations of the type $\Phi(x,\varphi(x),\varphi'(x),...,\varphi^{(n)}(x))=0$ when $\Phi$ is polynomial. (He defines a metatranscendental function as one which cannot satisfy an equation of this type for any $n$.) An interesting element is the use of the language of set theory and cardinalities: it is stated that for a fixed $n$ the set of solutions of such an equation is countable if the coefficients of the polynomial $\Phi$ are algebraic numbers but of cardinality continuum otherwise.\\

\end{enumerate}

\subsection{More on B\"ottcher's theorem} Introducing a new type of functional equation, called nowadays B\"ottcher equation, and solving it, is B\"ottcher's most important and most known achievement. The equation $Ff(z)=(F(z))^m$ and its generalizations can be considered under various assumptions about the given function $f$ (cf. \cite{KCG}, \cite{Re}). The case of a complex analytic $f$ satisfying $f(z)=\frac{z^m}{m!}f^{(m)}(0)+\frac{z^{m+1}}{(m+1)!}f^{(m+1)}(0)+...$ is interesting from the dynamical viewpoint since the solution $F$ provides an analytic conjugacy between $f$ and the monomial $z^m$ (and therefore it is often called B\"ottcher's coordinate). Because  only local behavior is considered, one can assume that the fixed point is at $0$. (By a change of coordinates, this in particular applies to  the important example of a polynomial $f$ of degree $m$ in a neighbourhood of its superattracting fixed point at infinity.) 
The following theorem, commonly known as B\"ottcher's theorem, summarizes  the result:
\begin{thm} Let $f(z) = a_mz^m+a_{m+1}z^{m+1}+..., \quad m \geq 2, a_m \neq 0$ be an analytic function in a neighborhood of $0$. Then there exists a conformal map $F$ of a neighborhood of $0$ onto the unit disk, $F(z)=z+bz^2+...$,  satisfying B\"ottcher's equation $Ff(z) = [F(z)]^m$.
\end{thm}

As we already mentioned, B\"ottcher proposed two (related) ways of constructing such a map (also called B\"ottcher's function): one using the quantity $\lim_{n \to \infty}\sqrt[m^n]{f_n(z)-z}$ and the other  using solutions to Gr\'evy's equation, which were obtained as path integrals. Both approaches were only sketched.  
Let us quote B\"ottcher's own summary of the calculations leading to his solution as given in  \S 34 of his Russian paper from 1904:\\

``It remains for us to study another important algorithm. We have the formula $B(z)=\frac{\lg (z_n-x)}{m^n}, \quad n=\infty$.  We write this formula in the  form $B(z)=\lg \sqrt[m^n]{(z_n-x)}$ and we get $N(z)=e^{B(z)}= \lim_{n \to \infty}\sqrt[m^n]{(z_n-x)}$. This algorithm satisfies the condition $Nf(z)=[N(z)]^m$, and since 
\[
B(z)=\lg (z-x)+ \frac{\lg \frac{f^{(m)}(x)}{m!}}{m-1}+\frac{A_0}{1}(z-x)+\frac{A_1}{2}(z-x)^2+...,
\]
then 
\[
N(z)=e^{B(z)}=e^{\lg (z-x)}\cdot e^{\frac{\lg \frac{f^{(m)}(x)}{m!}}{m-1}}\cdot e^{\frac{A_0}{1}(z-x)+\frac{A_1}{2}(z-x)^2+...},
\]
i.e., 
\[
N(z)=\sqrt[m-1]{f^{(m)}(x)/m!}\ (z-x)\biggl [1+\sum_{\nu=1}^{\infty}\frac{{\bigl(\frac{A_0}{1}(z-x)+\frac{A_1}{2}(z-x)^2+...\bigr)}^{\nu}}{\nu!}\biggr],
\]
from which it follows that 
\[
N(z)=\sqrt[m-1]{f^{(m)}(x)/m!}\ (z-x)+\frac{N^{(2)}(x)}{2!}(z-x)^2+\frac{N^{(3)}(x)}{3!}(z-x)^3+... .
\]

The region of convergence is the same as for the function $B(z)$, and if $z$ leaves the region, we compute $B(z)$ on the basis of the remark made in \S 31, and therefore we compute $N(z)$."\\

The details of the solution involving taking successive roots of iterates were supplied about 1920, independently by J.F. Ritt (\cite{Ri}) and P. Fatou (\cite{Fa}). In his article Ritt cites the Russian paper by B\"ottcher, while Fatou informs only that  ``the existence of this function (...) seems to have been first proved by Mr. B\"ottcher" (p. 189) and does not reference any of B\"ottcher's publications. (See \cite{Au} for more about early reception of B\"ottcher's theorem.) Here is Ritt's proof of B\"ottcher's theorem:\\

We can assume  that $a_m=1$ and consider a small disk in which  $f^p(x)$ approaches zero uniformly as p increases and such that $f$ has no zero other than $x = 0$. (In  this
disk, every iterate $f^p, \quad p = 1,2,...$, exists, and vanishes
only for $x = 0$, at which point it has a zero of order $m^p$.)  We can thus
select a sequence of functions $f^{1/m},f_2^{1/m^2},...,f_p^{1/m^p},...$ which is uniformly bounded. Now, $f_{p+q}(x)=[f^p(x)]^{m^q}(1+\epsilon_{p,q}(x))$, where $\epsilon_{p,q}(x)$ tends toward zero as $p$ increases, uniformly with respect to
$x$ and $q$. Consequently, $[f_{p+q}(x)]^{1/m^{p+q}}= [f^p(x)^{1/m^p}(1+\eta_{p,q}(x))$ where $\eta_{p,q}(x)$ approaches zero as $p$ increases, uniformly with respect to
$x$ and $q$. The functions in the selected sequence  are, for $p$ large enough, less than one 
in modulus, since they are roots of functions which approach zero, and it follows
 that they converge uniformly to an
analytic function $F$. Furthermore, since $[f_p(f(x))]^{1/m^p}=[[f_{p+1}(x)]^{1/m^{p+1}}]^m$, we find, as $p$  approaches infinity, that $Ff(x) = [F(x)]^m$.\\

 Fatou's proof is more detailed, but follows the same idea. It was reproduced in several references  mentioning B\"ottcher's  function, e.g. in \cite{Va}, \cite{Bl} and \cite{CG}. (The maps considered in \cite{Va} are holomorphic self-maps of the unit disk.) See also an account of this proof in \cite{AIR}.\\

In modern days new proofs of B\"ottcher's theorem were given. We will present several such proofs. \\

\textbf{J. Milnor's proof:} (\cite{Mi1}): This proof is in the classical spirit. The superattracting fixed point is placed at infinity rather that zero, i.e.,  the expansion of $f$ is 
\[
f(z) = a_mz^m+a_{m-1}z^{m-1}+...+a_0+a_{-1}z^{-1}+..., \quad |z|>R.
\]

By considering $z \mapsto \alpha f(z/\alpha)$ with $\alpha^{m-1}=a_m$ one can assume that $a_m=1$, so that $f(z)=z^m(1+\mathcal{O}(|1/z|)$ for large $z$. Substituting $e^z =Z$ lifts $f$ to a continuous map $F(Z)= \log f(e^Z)$, uniquely defined for $\Re(z) > \log R$ up to addition of some multiple of $2\pi i$. The lift $F$ can be chosen so that $|F(Z)-mZ|<1$ in some half-plane $\Re(z) > \sigma$ and this half-plane is mapped to itself by $F$. If $Z_0 \mapsto Z_1 \mapsto ...$ is any orbit under $F$ in this half-plane, then $|Z_{k+1}-mZ_k|<1$ and $|W_{k+1}-mW_k| < 1/m^{k+1}$, where $W_k := Z_k/m^k$.  Hence the sequence of holomorphic functions $W_k=W_k(Z_0)$ converges uniformly as $k \to \infty$ to a holomorphic function $\Phi$. This function  satisfies $\Phi(F(z))=m\Phi(Z)$ and $\Phi(Z+2\pi i)=\Phi(Z)+2\pi i$. Therefore the mapping $\phi(z) = \exp(\Phi(\log z))$ is well defined near infinity and satisfies $\phi(f(z))=(\phi(z))^m$.\\

Here is a non-classical proof, which  uses  quasiconformal maps to extend a partially defined topological conjugacy and improve its regularity. \\

\textbf{M. Lyubich's proof:} (\cite{L}) Consider a small disk $U_\varepsilon$ and let $W$ denote the connected component of $f^{-1}U_\varepsilon$ containing $0$. Let $V=\overline{W\setminus U_\varepsilon}$ and let $r \in (0,1)$. Consider a diffeomorphism $h_0: V \mapsto A[r^k,r]$ such that $h_0(f(z))=h_0(z)^m, \quad z \in \partial W$. This diffeomorphism extends to a homeomorphism $h: W \mapsto U$ for which $h(f(z))=h(z)^m, \quad z \in W$. The conformal structure $h_*\sigma$ extends naturally  to a conformal structure on the unit disk $U$ which is invariant under the map $G: z \mapsto  z^m$ in a neighborhood of an arbitrary point different from zero. By the measurable Riemann mapping theorem, there exists a quasiconformal homeomorphism $\psi:\bar{ U} \mapsto \bar{U}$ with $\psi(0)=0, \psi(1)=1$ and $\psi_*\mu=\sigma$. Then the mapping 
 $g=\psi \circ G \circ \psi^{-1}$ is locally conformal outside $0$ and is an $m$-sheeted covering of $U^*$ by itself. Then $g: z \mapsto z^m$.\\

We mention that there is  another proof of B\"ottcher's theorem using quasiconformal maps (combined with holomorphic motions), due to Y. Jiang (\cite{Ji}).\\

Yet another proof uses operator theory after observing that the exponential of the B\"ottcher function is an eigenvector of a composition operator:\\ 

\textbf{T. Gamelin's proof:} (\cite{Ga}) Let $U_\delta$ be a disk centered at $0$ with radius $\delta$ such that $0 < \delta < 1/2$. Let $A_0$ be the space of functions $g$ analytic in $U_\delta$ satisfying $g(0)=0$ and let $T_0$ be the restriction to $A_0$ of the composition operator $(Tg)(z)=g(f(z))$.  We can assume that $|f(z)| \leq 2|z|^m$ in $U_\delta$, so $|f^n(z)| \leq (2\delta)^{m^n}$. If $g \in A_0$ satisfies $\|g\|\leq 1$, then by Schwarz lemma $|(T_0^ng)(z)| \leq |f^n(z)| \leq (2\delta)^{m^n}$, so $\|T_0^n\| \leq (2\delta)^{m^n}$. Hence the Neumann series $(\lambda I -T_0)^{-1} =\sum_{k=0}^{\infty} \frac{T_0^k}{\lambda^{k+1}}$ converges for all $\lambda  \neq 0$. Note that B\"ottcher equation is equivalent to the resolvent equation $(mI-T_0)G=h$, where $G(z)=\log(F(z)/z)$, $h(z)=\log(f(z)/z^m)$. It is thus enough to solve the resolvent equation for $G \in A_0$ and then set $F(z)=z\exp(G(z))$. Such a solution can be obtained using Banach contraction principle.\\

The operator-theoretic approach was also applied in \cite{Kl}, where the (somewhat more general) B\"ottcher function was obtained as the unique fixed point of a certain contraction operator. For the relation with composition operators, see \cite{Ne}.\\

A construction of the inverse map to B\"ottcher's function for the family of quadratic maps $f_\lambda: z \mapsto (z-\lambda)^2$ by an iterational process based on Caratheodory's theorem on convergence of domains can be found in \cite{BGH}.

It should be mentioned that the classical method of proving B\"ottcher's theorem was used recently in \cite{Fav} to give the classification (up to holomorphic conjugacy) of attracting rigid germs in $\mathbb{C}2$ (which was subsequently applied in \cite{Ru}), and  in \cite{KR} in the study of local stable manifolds of a dominant meromorphic self-map $f: X \dashrightarrow X$, where $X$ is a compact K\"ahler manifold of dimension $n > 1$.   As far as other applications and generalizations are concerned, see e.g. \cite{Pog}. An analog of B\"ottcher's theorem for a class of transcendental entire functions was proved in \cite{Re}, while an analog of B\"ottcher coordinate for mappings of the form $f= g \circ h$, where $g$ is a polynomial of degree $d \geq 2$ and $h$ is an affine mapping of the complex plane to itself, was found in \cite{FF}. 
The existence of a local B\"ottcher coordinate for a holomorphic map in several complex variables with a superattracting fixed point was proved in \cite{BEK}; some special cases and examples were introduced and studied in \cite{U}. An analog of B\"ottcher coordinate for a map with a superattracting hyperplane was considered in \cite{BJ}.\\

\section{Appendix: archival materials}

The documents reproduced here come from: L. B\"ottcher, Teka osobowa (personal file), fond 27, op.4-40, Lvov District Archive, and 

L. B\"ottcher, Universit\"atsarchiv, Lepizig, Phil.Fak., Prom., 714, Bl.  7. We thank these institutions. All translations are by the second author.

\subsection{Programmes in mathematics in the c.k. Polytechnic School in Lvov during L. B\"ottcher's studies}

In the Division of Machine Construction of the Polytechnic School B\"ottcher attended the following lectures in mathematics: Mathematics course I, 6 hours of lectures per week in the winter and summer semester (dr Placyd Dziwi\'nski); descriptive geometry, 5 hours of lectures per week in the winter and summer semester, 10 hours of repetitory in the winter semester (dr Mieczys\l aw \L azarski);  repetitory in elementary mathematics, 2 hours per week in the winter semester and 10 hours per week in the summer semester (dr Placyd Dziwi\'nski); mathematics course II, 6 hours per week in the spring and summer semester (dr W\l adys\l aw Zaj\c aczkowski); repetitory in higher mathematics, 2 hours per week in the winter and summer semester (dr W\l adys\l aw  Zaj\c aczkowski); drawings in descriptive geometry, 10 hours of exercises in the winter and summer semester.\\

The programme of course I was the following: {\it Foundations of higher analysis}- an introduction to analysis: theory of operations; infinite series and products; algebraic equations; determinants and methods of elimination; variable quantities and their functions; differential calculus: differentials and derivatives of functions of one and several variables; Taylor and Maclaurin formulas; indeterminate symbols; maxima and minima; tangency and curvature of planar and spatial curves; tangency of surfaces; integral calculus: methods of integration; integrals of rational, algebraic and transcendental functions; approximate methods of computing integrals; multiple integrals; rectification and quadrature of curves; quadrature and cubature of surfaces.\\
{\it Analytic geometry:} systems of coordinates in the plane and in space; formulas of planar and spherical trigonometry;  points, straight lines and planes;  geometrical loci; homogeneous coordinates; the cross-ratio and involution; theory of curves and surfaces of second degree.\\

The programme of course II was the following: {\it Higher analysis}  -- the theory of definite integrals; methods of computing definite integrals; multiple definite integrals; Euler integrals; Fourier integrals and series; theory of functions of an imaginary variable; differentials and integrals of complex  functions; general properties of analytic functions; theory of differential equations; formation of differential equations; the  theory of Jacobians;  integration of ordinary differential equations of the first order and higher orders, especially linear ones; integration of a system of ordinary differential equations; integration of partial differential equations of the first order, linear and general ones with three variables; foundations of the calculus of variations.
{\it General theory of curves and surfaces} -  osculation and curvature of non-planar curves  and surfaces; ruled surfaces; quadric surfaces; curves on surfaces; curvature, geodesic and asymptotic lines; cubature and quadrature of a surface. \\

The programme in descriptive geometry was the following:
{\it Methods of descriptive geometry} -- central projections; series of points and pencil of rays; unicursality  of pencils and series;  theory of curves of second degree; collineation, similarity,  involution, congruence and  symmetry of planar systems; collineation and similarity of spatial systems; orthogonal projections;  axonometry;
{\it Theory of curves and surfaces in general}- non-planar curves and expandable   surfaces; cones and cylinders, non-planar curves of the 3-rd and 4-th order;  the helix and the expandable helical surface. 
{\it Theory of oblique surfaces} - the hyperboloid with one  sheet; the hyperbolical paraboloid; helical  surfaces;
{\it Theory of non-ruled surfaces of the 2-nd order} - the sphere; rotational surfaces of the second order as collineates of the sphere; Tri-axial surfaces of second order as objects related with rotational surfaces of second order;
{\it Theories of rotational surfaces and envelopes} - constructions of proper and cast shadows and lines of equal illumination on a surface. \\

There were also repetitory courses in elementary and higher mathematics. They were respectively recitations corresponding to course I and  supplementary chapters in course II.\\ 

\newpage

\subsection{B\"ottcher's mathematical studies in Leipzig}


\begin{center}

Matriculation in the university courses at Leipzig 
\includegraphics[height=17cm]{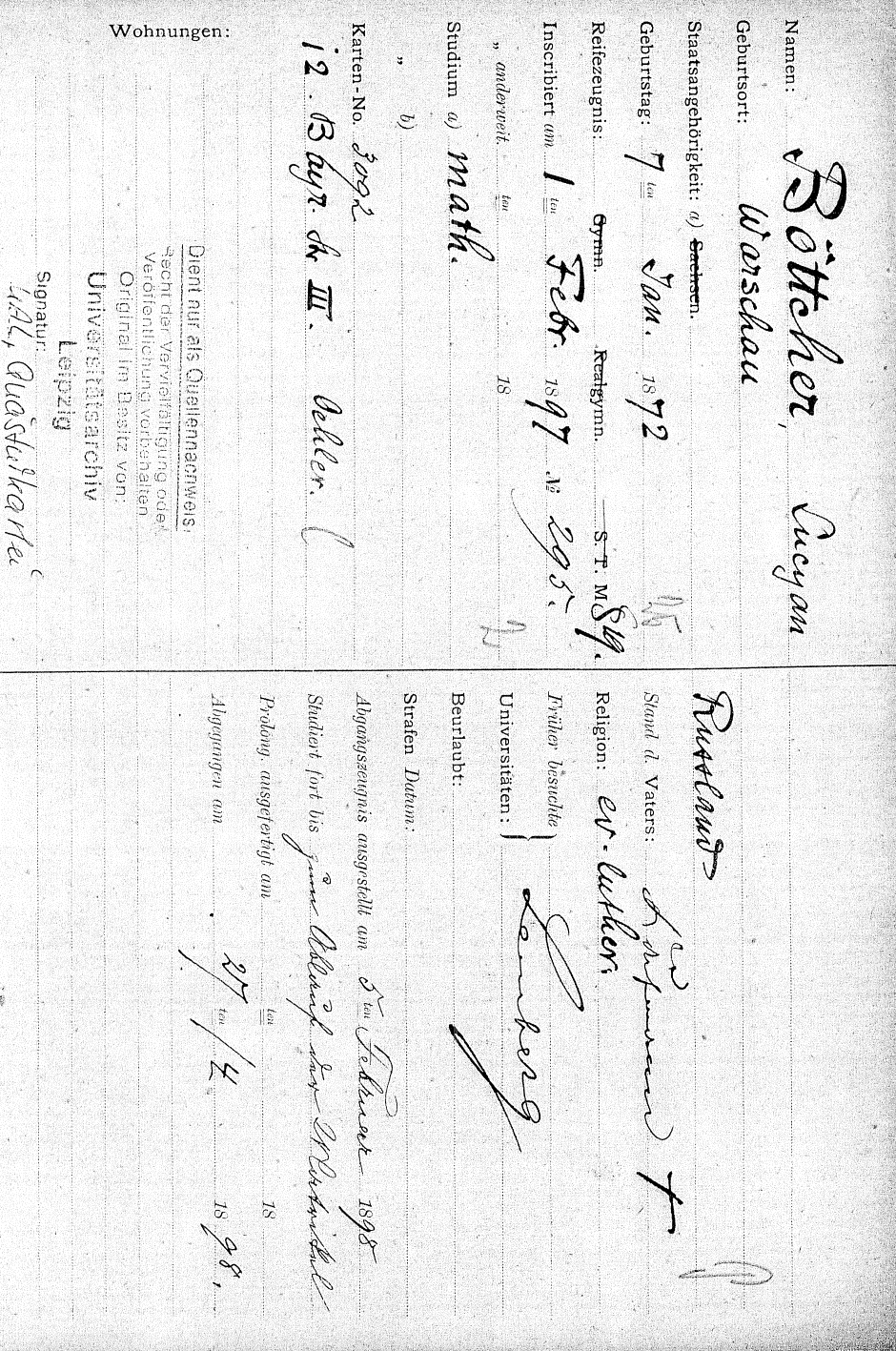}

\end{center}

 

\begin{center}

\newpage

L. B\"ottcher's CV in Latin, printed with his thesis

\includegraphics[height=15cm]{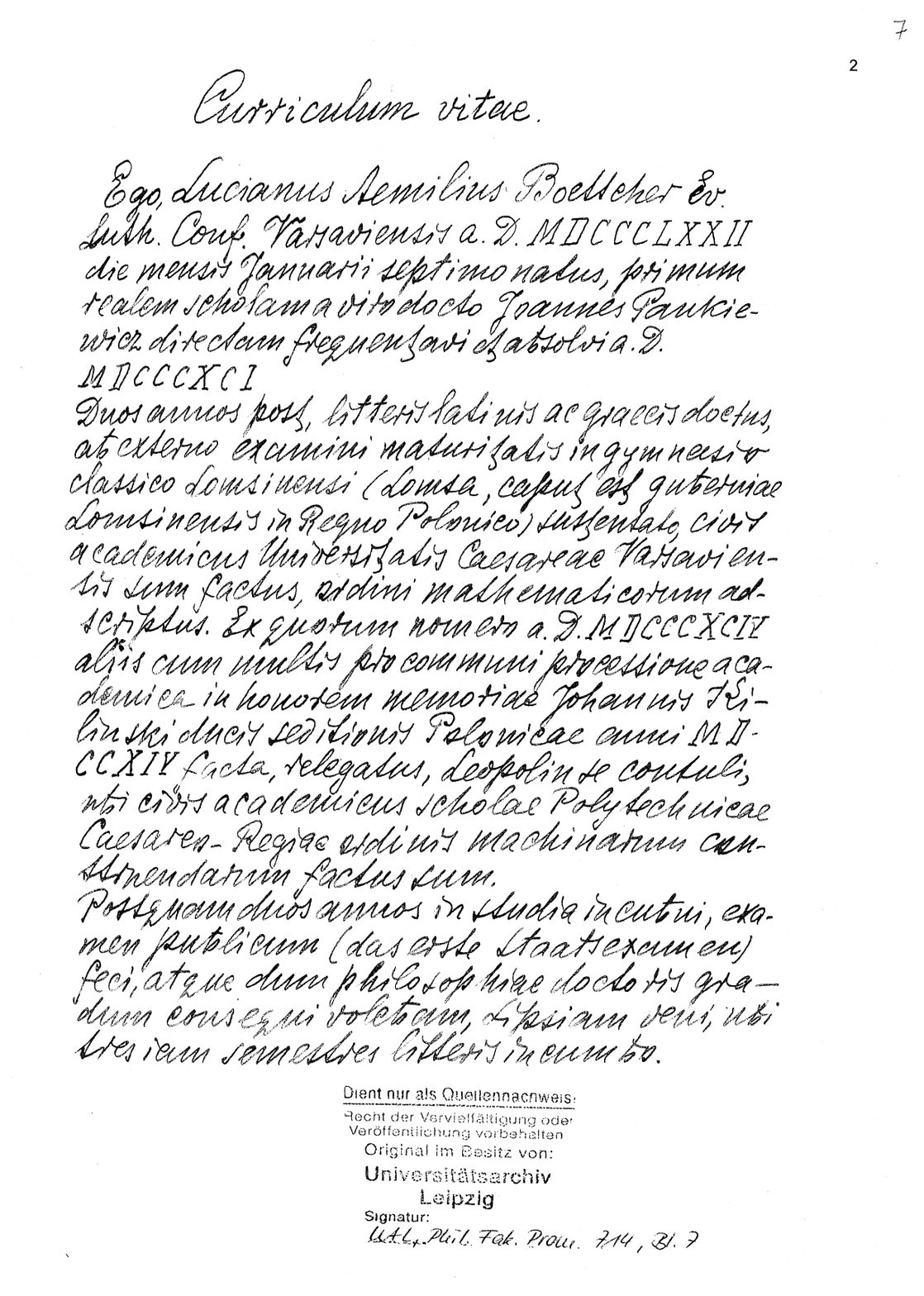}

\end{center}



\newpage

Formal evaluation of B\"ottcher's thesis and examinations

\vspace{1cm}
\begin{center}

\includegraphics[height=15cm]{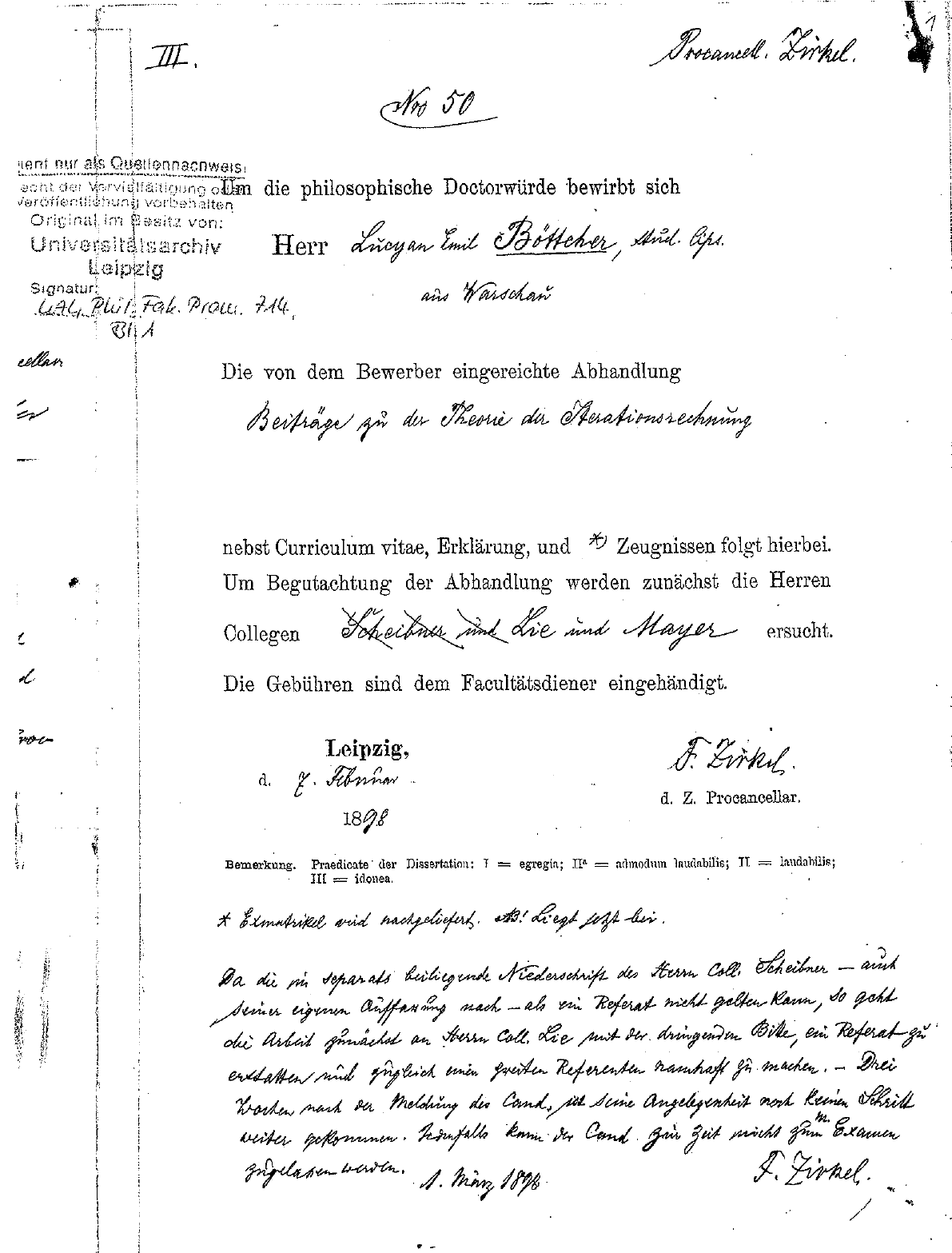}

\end{center}



\newpage

\begin{center}

Sophus Lie's report on B\"ottcher's thesis

\includegraphics[height=15cm]{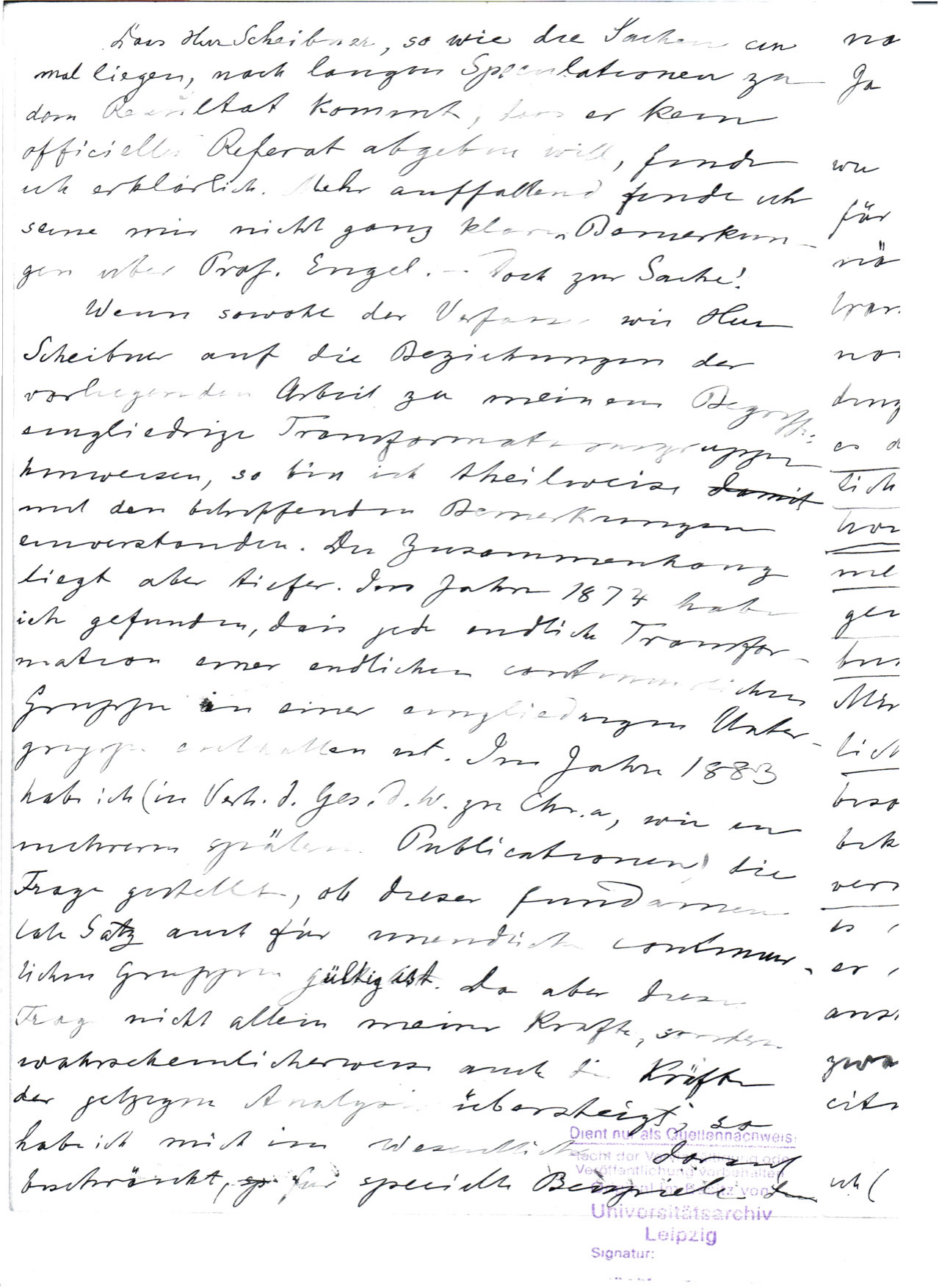}

\includegraphics[height=15cm]{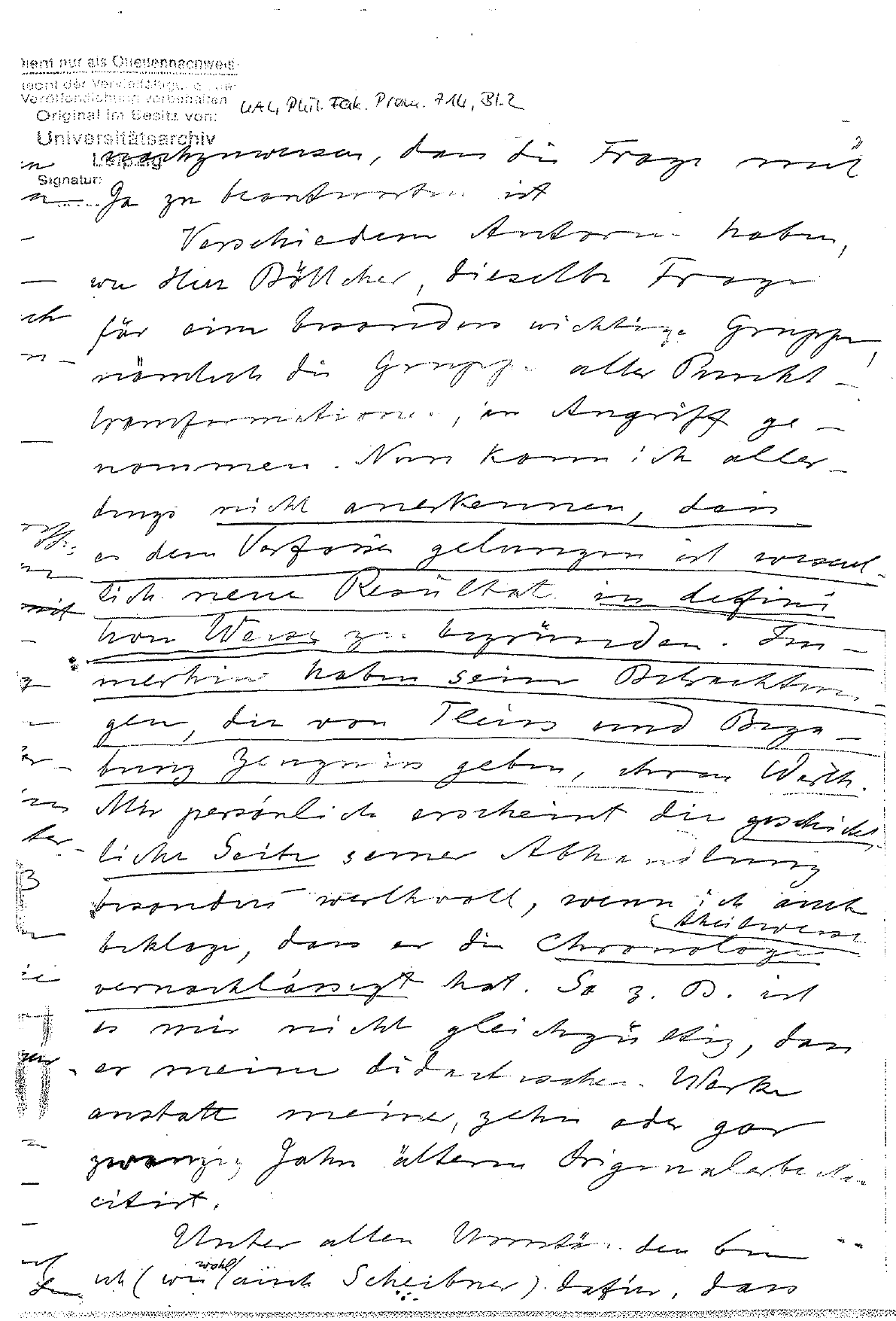}

\includegraphics[height=15cm]{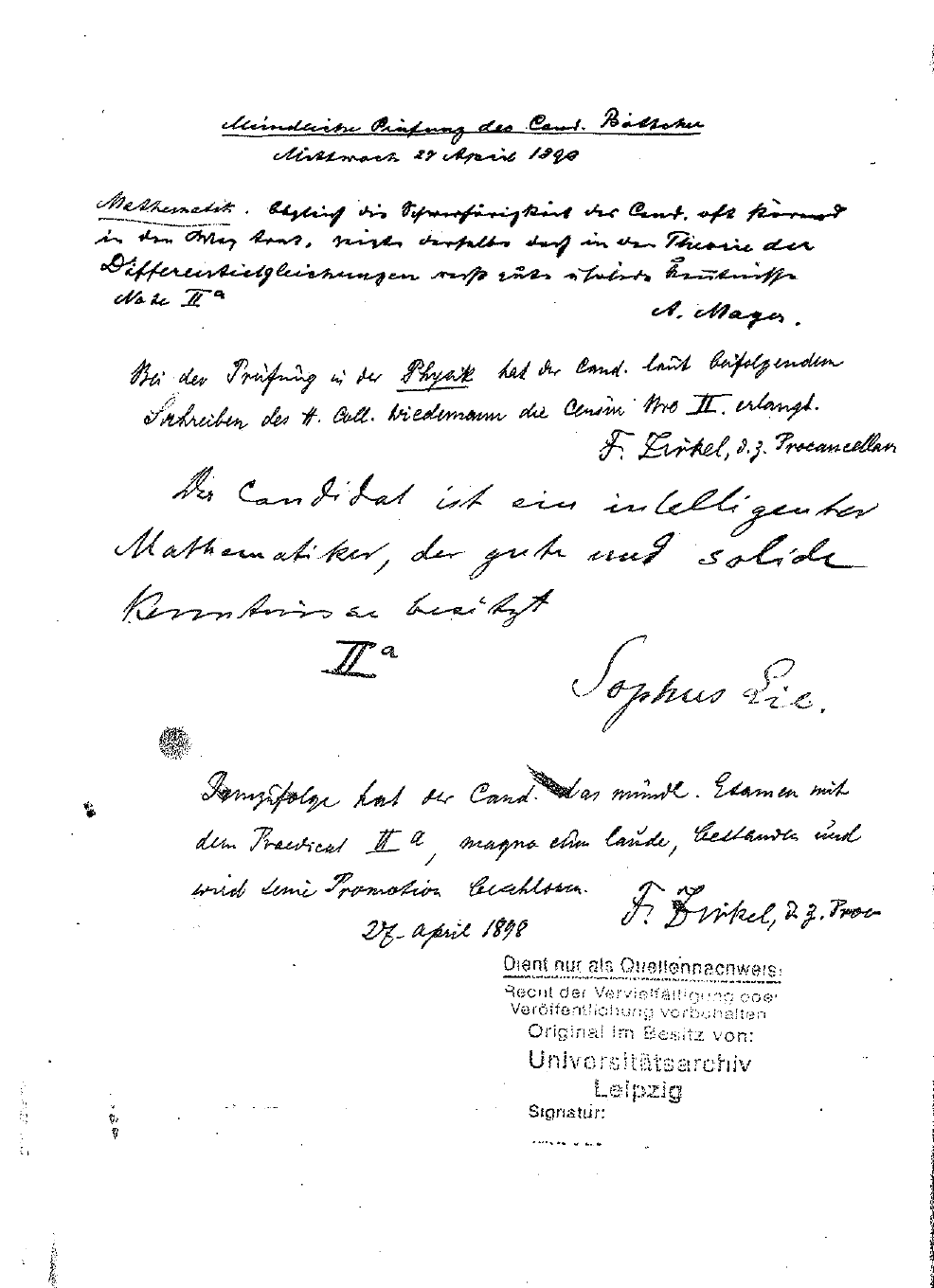}

\end{center}


\newpage

\subsection{L. B\"ottcher's habilitation at the Lvov university}

Attempt 1, 1901\\

Sources: Lvov District Archive\\
Lvov Polytechnic School\\
L. B\"ottcher's personal file\\
27.4.40\\


B\"ottcher's CV, with information about his studies in Warsaw and Leipzig and his scientific activities afterwards

\begin{center}

\includegraphics[height=15cm]{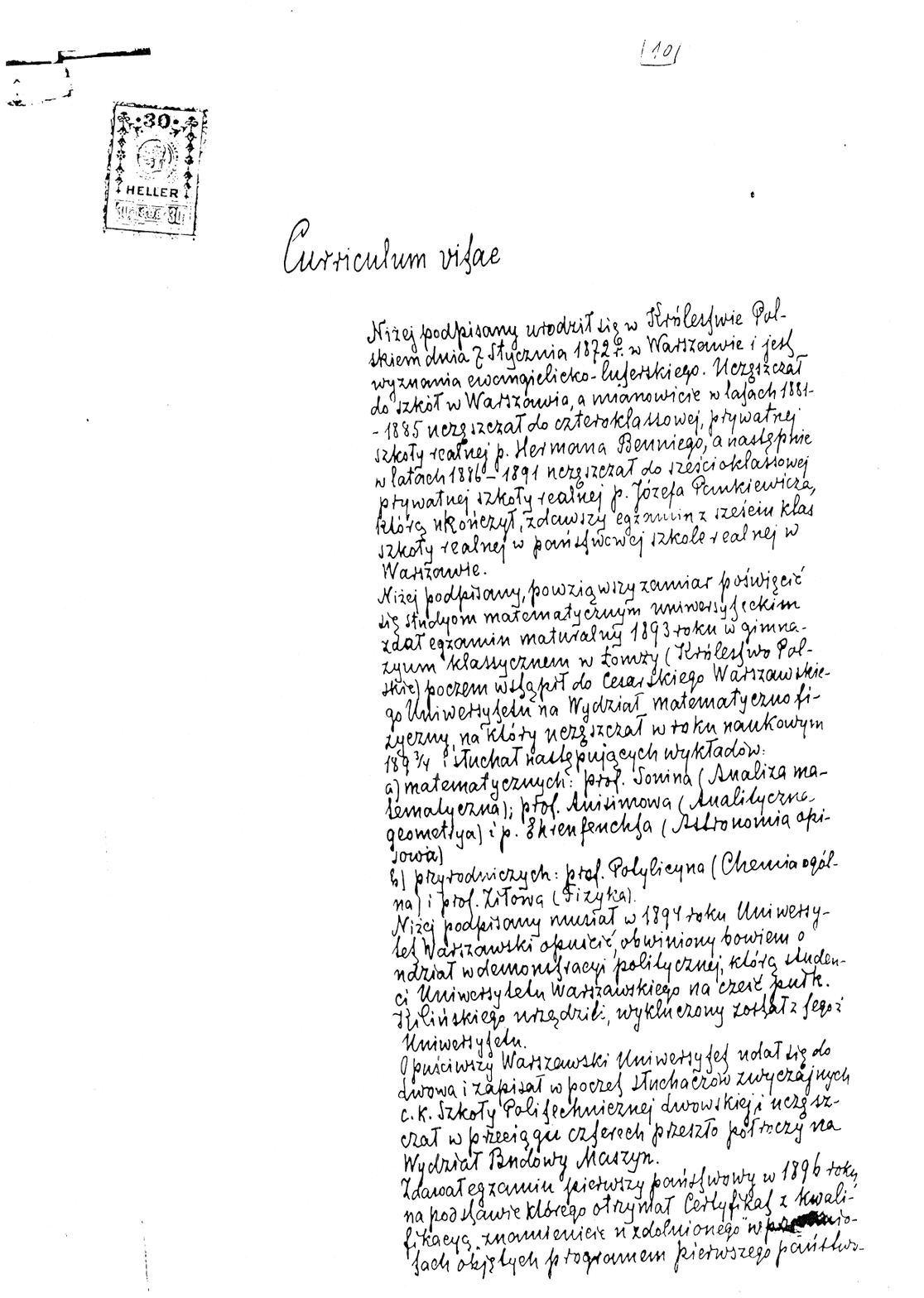}

\end{center}
  

\begin{center}


\includegraphics[height=15cm]{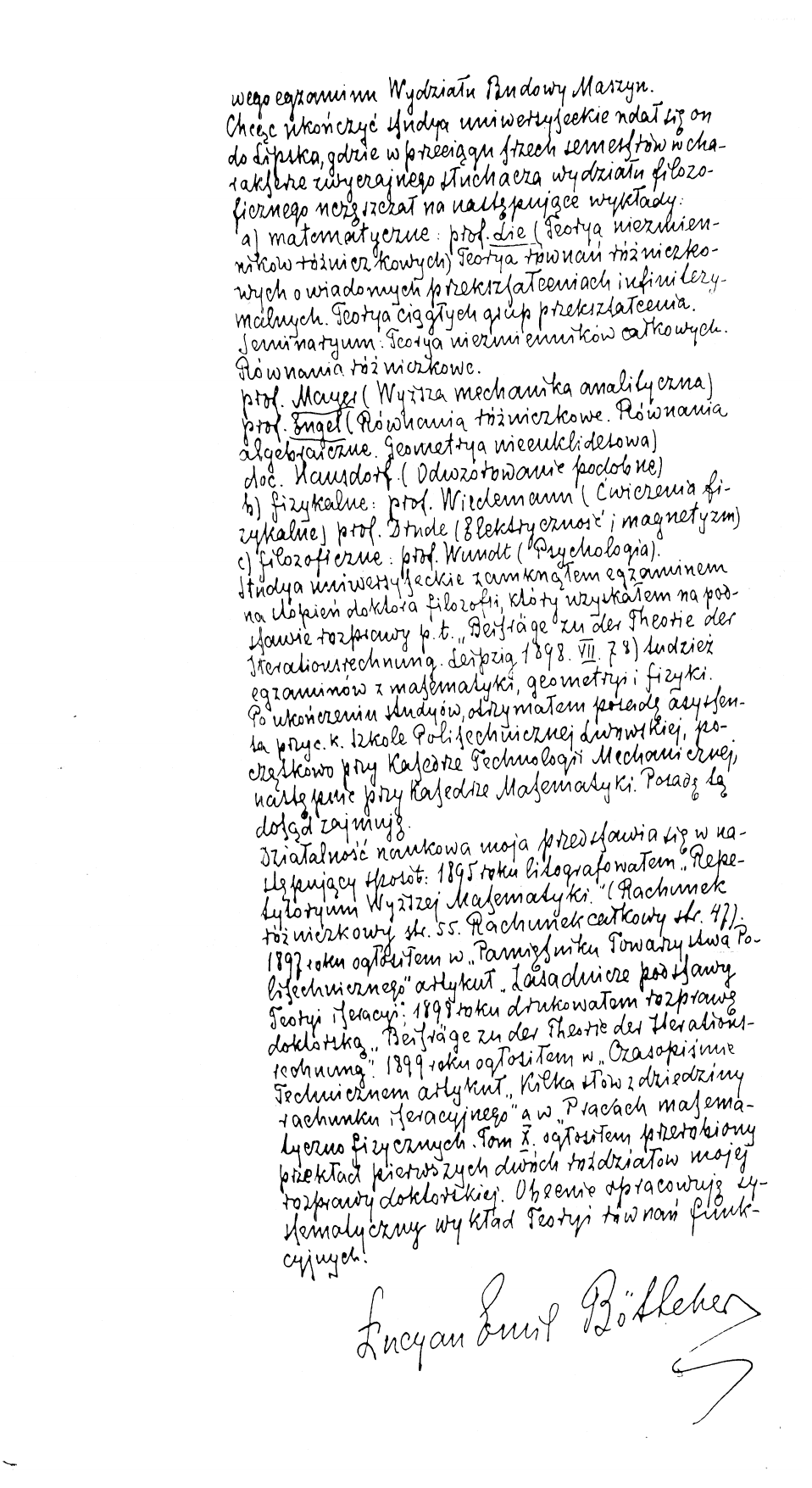}

\newpage
\end{center}

B\"ottcher's applications for admission to veniam legendi and habilitation
\begin{center}

\includegraphics[height=15cm]{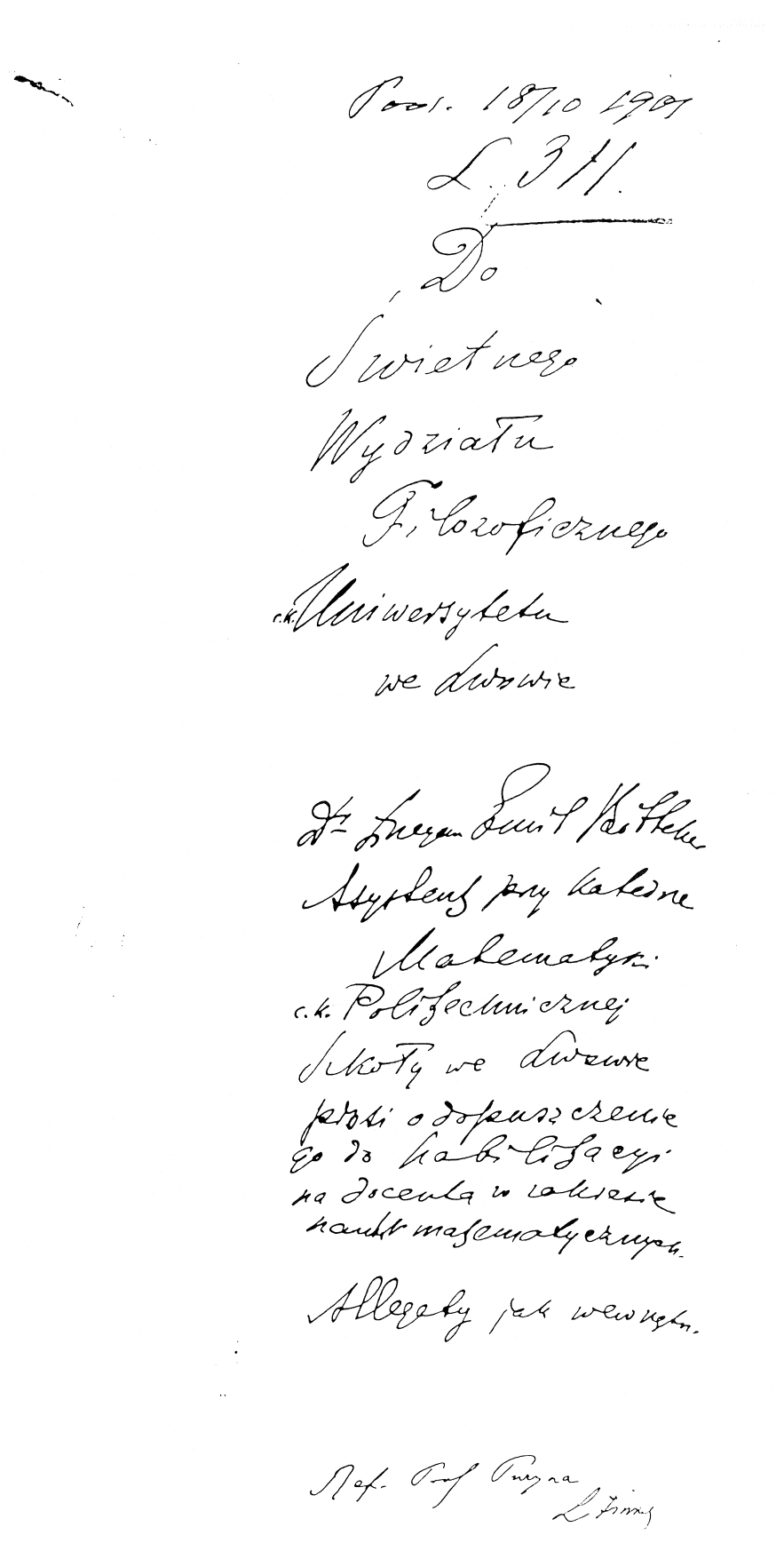}

\end{center}

\begin{center}


\includegraphics[height=15cm]{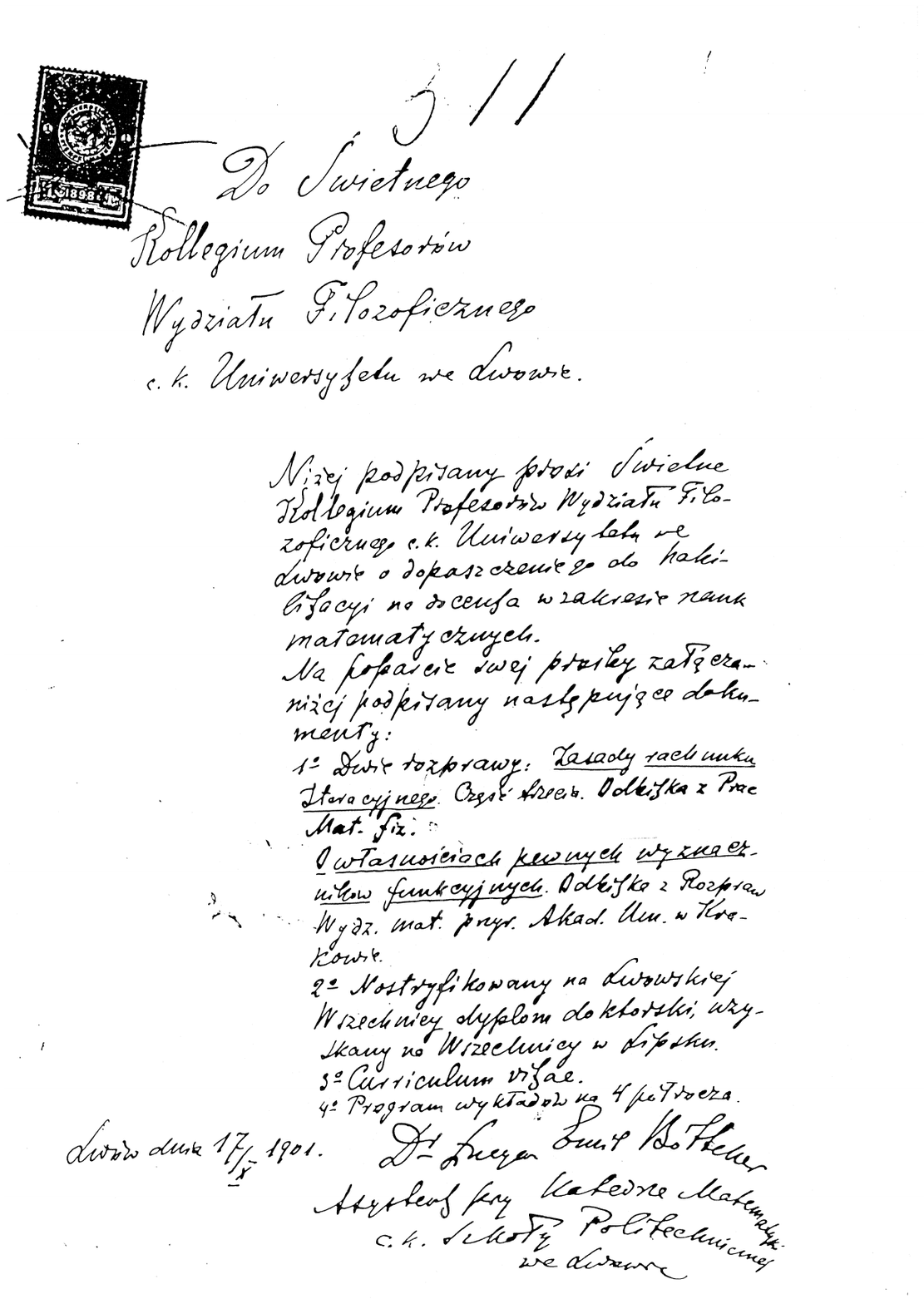}

\end{center}


To the Honourable Collegium of Professors\\
of Philosophical Faculty\\
of the c.k. University in Lvov\\

The undersigned requests from the Honourable Collegium of Professors
of the Philosophical Faculty 
of the c.k. University in Lvov  to be admitted to habilitation as a {\it docent} in mathematical sciences. In support of his request, the undersigned encloses the following documents:\\

1. Two papers: {\it Principles of calculus of iterations. Part III Applications of the theory of convergence of iterations to solving elementary functional equations}, an offprint from {\it Prace mat.--fiz.}\\
{\it On properties of some functional determinants}, an offprint from {\it Rozprawy Wydz. Mat. Przyr. Akademii Umiej\c etnosci w Krakowie}\\
2. Doctoral diploma, nostrificated at the Lvov University, obtained at the Leipzig University\\
3. Curriculum vitae\\
4. A programme of lectures for 4 semesters.\\
Lvov 17 X 1901\\
Signed Dr. \L ucjan Emil B\"ottcher\\
Assistant in the Chair of Mathematics of the Polytechnic School in Lvov.\\
Below, the Dean Finkel directs a letter to Prof. Puzyna\\


To the Honourable Faculty of Philosophy of the  University in Lvov\\
Dr.  \L ucjan Emil B\"ottcher Assistant in the Chair of Mathematics of the Polytechnic School in Lvov requests to be admitted to habilitation as a {\it docent} in mathematical sciences.\\

 
I invite Honourable Colleagues, Members of the Committee for Dr. B\'ottcher's habilitation, [to the meeting] which will take place tomorrow, Monday January 13, at the hour 4. in the afternoon. Lvov, 12 I 1902.\\
Signed L. Finkel\\
Received:\\
HH. Prof. Puzyna\\
HH. Prof. Rajewski\\
HH. Prof. Smoluchowski\\

The meeting on 13 I 1902\\
Present: the Dean in charge\\
Prof. Puzyna, Smoluchowski, Rajewski.\\

\textbf{The motion:} So far, the committee decided to wait for publication of the paper {\it Principles of the iterational calculus, Part III}.\\
Professors' signatures.\\


\newpage


\begin{center}

Unanimous motion of the committee regarding Dr. B\"ottcher's habilitation

\includegraphics[height=15cm]{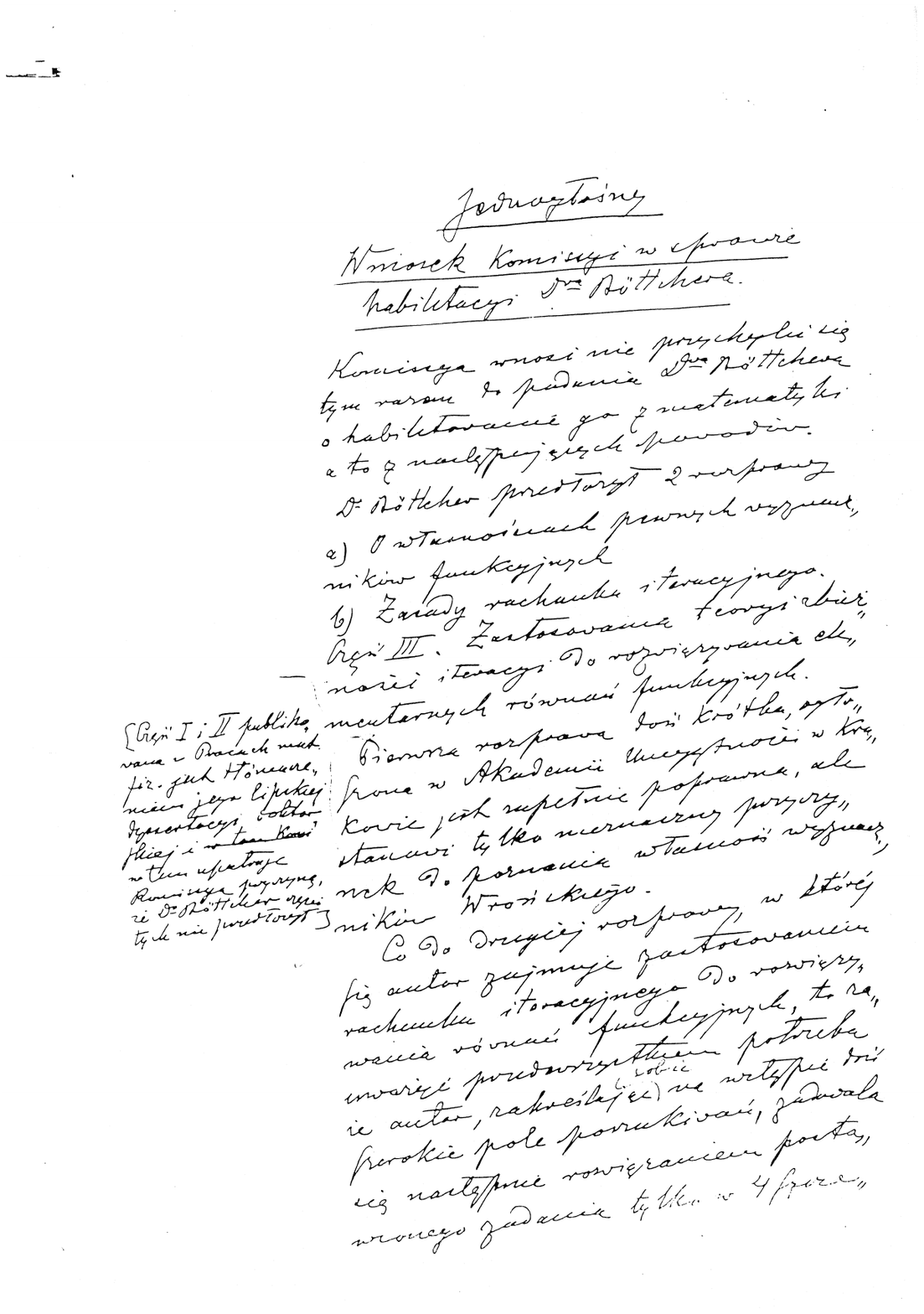}

\includegraphics[height=15cm]{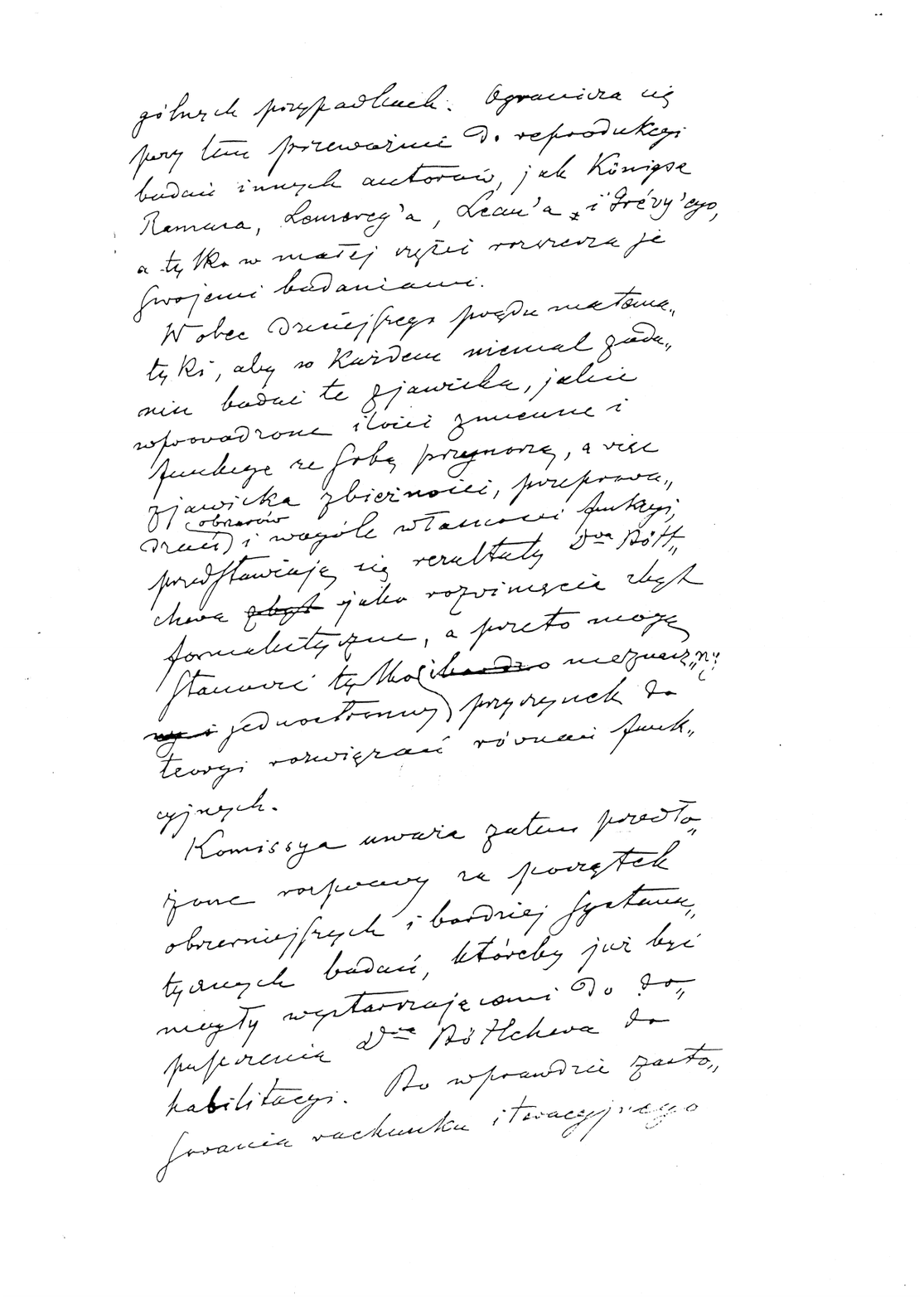}

\includegraphics[height=15cm]{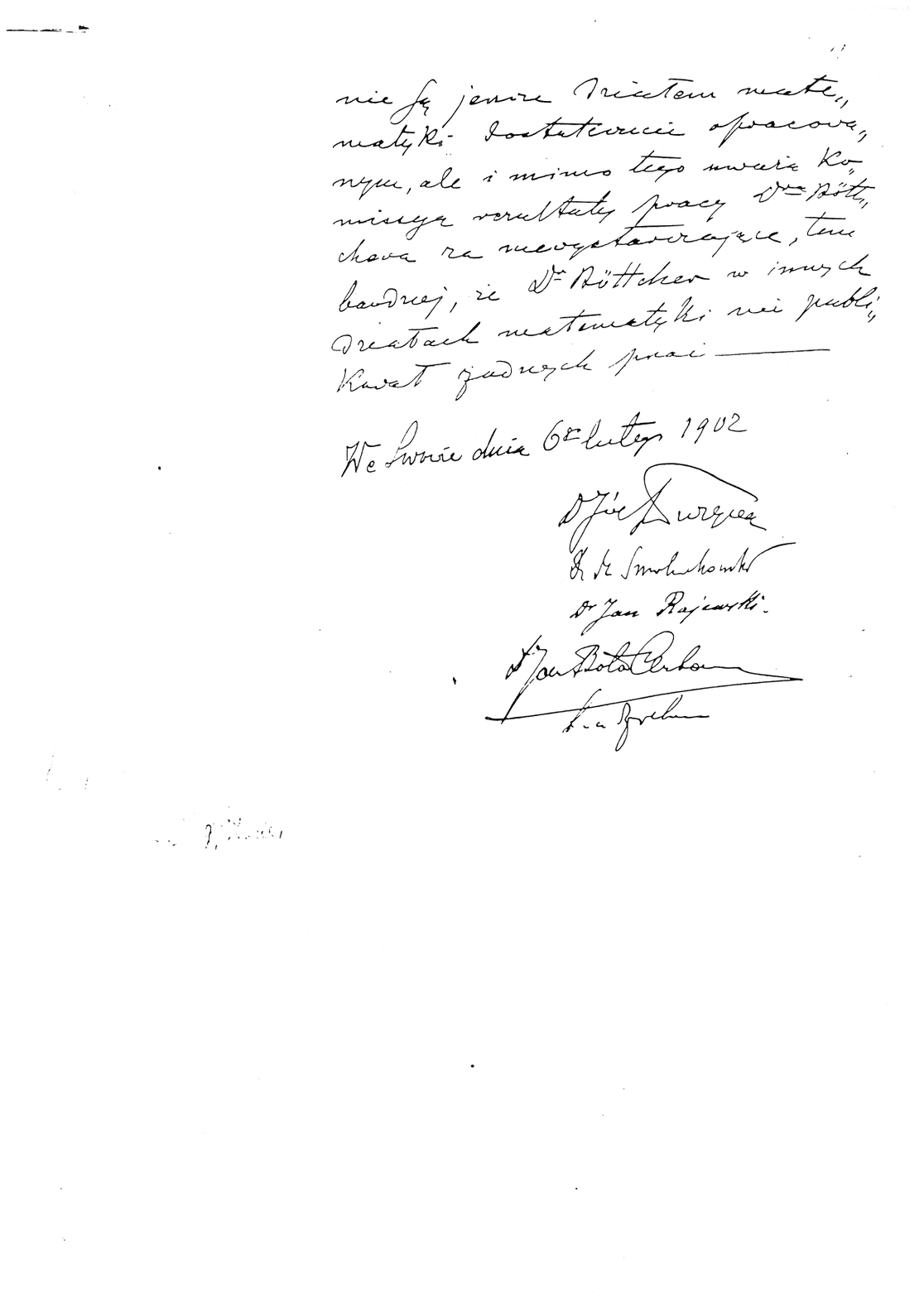}

\end{center}


The committee motions not to approve this time  Dr. B\"ottcher's application for habilitation in mathematics, because of the following reasons.\\
Dr. B\"ottcher submitted two papers\\
a) On properties of some functional determinants\\
b) Principles of calculus of iterations. Part III Applications of the theory of convergence of iterations to solving elementary functional equations.\\

[The margin note: parts I and II published in Prace mat.--fiz.[sic] are a translation of his Leipzig dissertation, and the Committee views it as a reason why Dr. B\"ottcher did not submit these parts.]\\

The first paper, fairly short, published by the Academy of Skills in Krak\'ow, is completely correct, but is only a small contribution to the study of properties of Wro\'nskian determinants.  As regards the second paper, in which the author deals with applications of  calculus of iterations to solving functional equations, one needs first of all to note that the author, outlining at the beginning quite a wide range for his investigation, is then satisfied with solving the posed problem only in 4 special cases. In this, he restricts himself mainly to reproduction of considerations of other authors, such as K\"onigs, Ramus, L\'emeray, L\'eau, and Gr\'evy, extending these considerations only a little by his own investigations. In view of the current trend of mathematics to study in almost all problems the phenomena brought by introduction of variable quantities and functions, that is, the phenomena of convergence, transformations of domains, and properties of functions in general, the results of Dr. B\"ottcher seem to be too formalistic developments, and therefore can be only one-sided contributions to the theory of solutions of functional equations.\\
  Thus the Committee views the submitted papers as a beginning of broader and more systematic investigations which could become sufficient for admitting Dr. B\"ottcher to habilitation. While applications of calculus of iterations are not yet a well developed area of mathematics, the Committee nevertheless regards Dr.  B\"ottcher's results as insufficient, even more so because Dr. B\"ottcher did not publish any works in other areas of mathematics.\\
  Lvov, February 6, 1902\\
  Signed:\\
  Dr. J\'ozef Puzyna\\
  Dr. M. Smoluchowski\\
  Dr. Jan Rajewski\\
\newpage

Attempt 2, 1911\\


\begin{center}


\includegraphics[height=15cm]{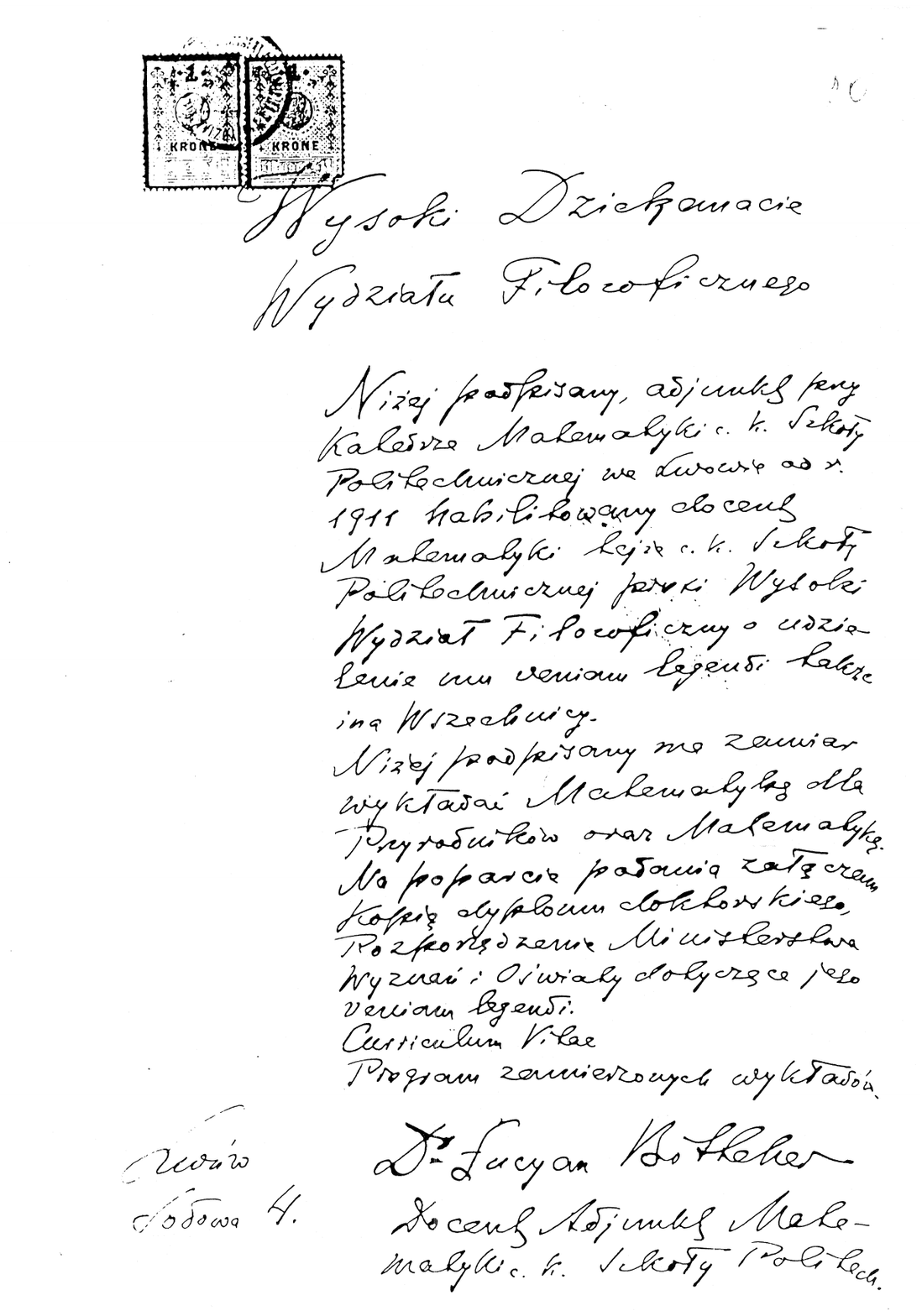}

\end{center}


  Honourable Dean's Office\\
  of the Philosophical Faculty\\
  The undersigned, an {\it adiunkt} in the Chair of Mathematics c.k. Polytechnic School in Lvov, a habilitated {\it Docent} of the same c.k. Polytechnic School since 1911, requests that the Honourable Philosophical Faculty grant me veniam legendi also at the University.\\
  The undersigned  intends to conduct lectures in Mathematics for Naturalists and Mathematics.\\
  To support the application I enclose a copy of the doctoral diploma, the decree of Ministry of Denominations and Enlightment regarding its veniam  legendi.\\
  Curriculum Vitae\\
  The programme of planned lectures.\\
  Lvov, Sodowa 4, Dr. Lucjan B\"ottcher, Docent Adiunkt of Mathematics of c.k. Polytechnic School\\



  
  Attempt 3, 1918\\
  Dr. Lucjan B\"ottcher\\
  Application for veniam legendi in mathematics sent in the meeting on ... to the elected Committee\\
  

\begin{center}


\includegraphics[height=15cm]{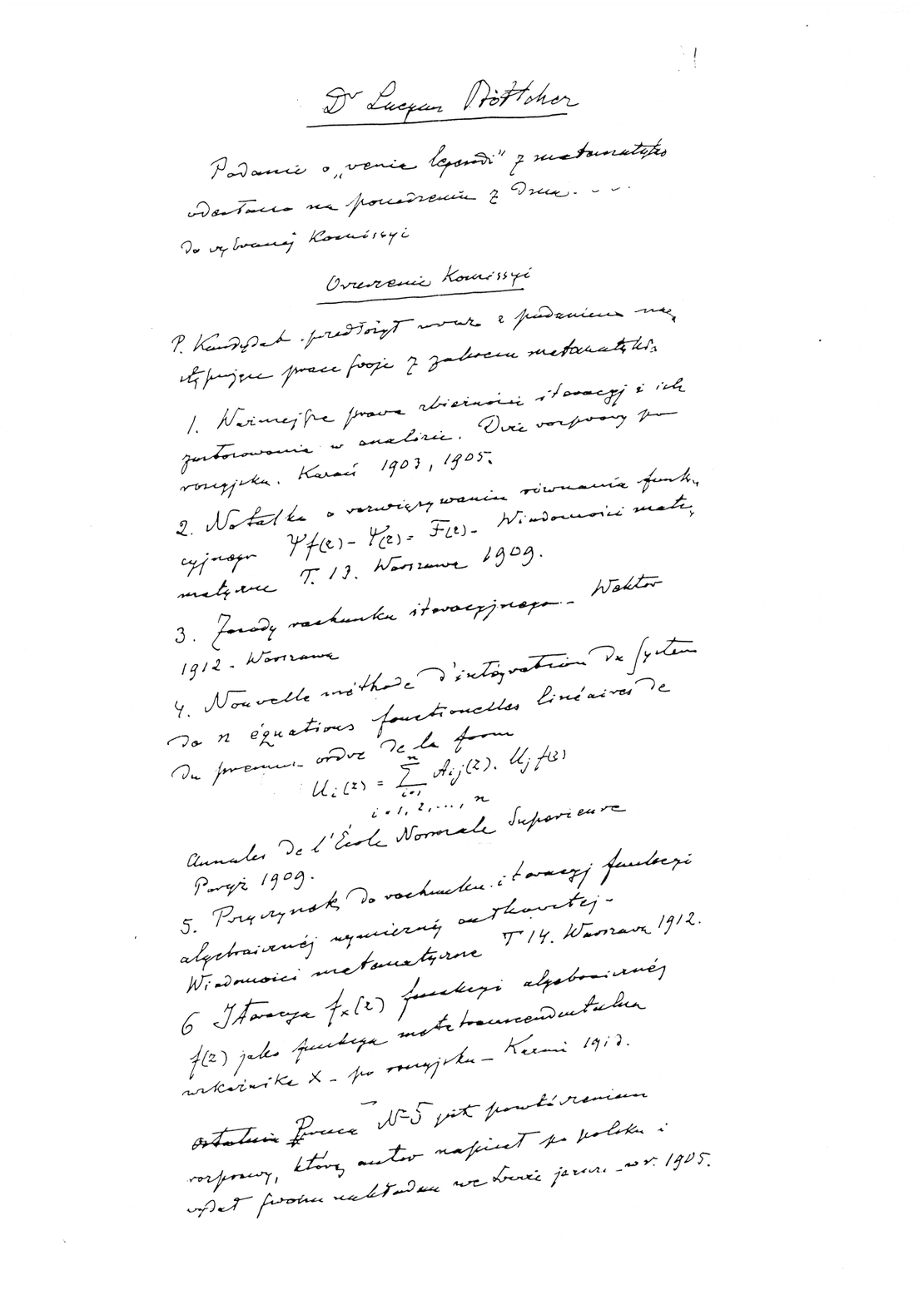}

\end{center}



\begin{center}

The Committee's Decision

\includegraphics[height=15cm]{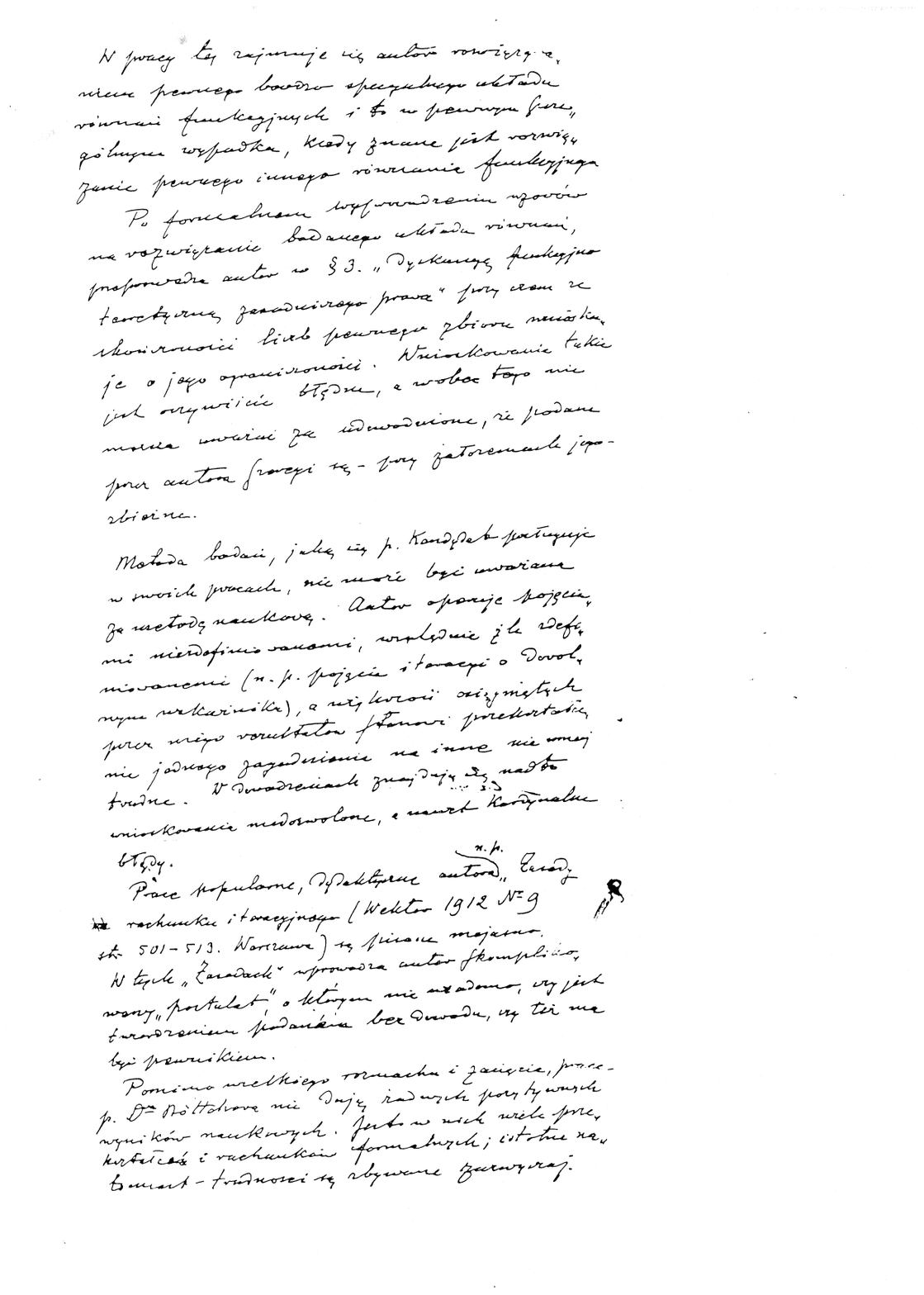}

\end{center}

 

  The Candidate submitted along with the application the following works in mathematics:\\
  1) Major laws of convergence of iterations and their applications in analysis. Two papers in Russian, Kazan, 1903, 1905.\\
  2)  A note of solving the functional equation $\Psi f(z)-\Psi(z)=F(z)$, {\it Wiadomo\'sci Matematyczne}, vol. 13, Warsaw 1909\\
  3) Principles of iterational calculus, {\it Wektor}, 1912, Warsaw\\
  4) Nouvelle m\'ethode d'int\'egration d'un syst\`eme de $n$ \'equations fonctionelles lineair\'es du premier ordre de la forme $U_i(z)=\sum_{j=1}^{j=n}A_{i,j}(z)U_jF(z)$, Annales l'Ecole Normale Sup\'erieure, Paris, 1909\\
  5) A contribution to the calculus of iteration of a rational entire function, {\it Wiadomo\'sci Matematyczne}, vol.14, Warsaw 1912\\
  6) Iteration $f_x(z)$ of an algebraic function $f(z)$ metatranscendental in the index $x$, in Russian, Kazan 1912\\
  
  The paper no. 5 duplicates one written by the author in Polish and self-published already in 1905. In it, the author deals with solving a very specific system of functional equations, in a very particular case in which a solution of some other functional equation is known. After formal deduction of formulas for a solution to the system of equations under investigation the author proceeds to give in \S  3 "A functional-theoretical discussion of the fundamental law", concluding boundedness of a certain set from finiteness of the numbers in it. Such reasoning is obviously erroneoeus, and therefore one cannot consider it to be proven that the series given by the author are-- under his conditions--convergent.\\

The method used by the Candidate in his works cannot be considered scientific. The author works with undefined, or ill-defined, notions (e.g., the notion of an iteration with an arbitrary exponent), and the majority of the results he achieves are transformations of one problem into another, no less difficult. In the proofs there are moreover illegitimate conclusions, or even fundamental mistakes. The author's popular, instructional works, e.g. "Principles of iterational calculus" ({\it Wektor} 1912, no. 9, pp.  501-513, Warsaw), are written in an unclear manner. In these "Principles" the author introduces a complicated new "postulate", about which it is not known whether it is a theorem stated without a proof or it is supposed to be an axiom.\\

Despite great verve and determination, Dr. B\"ottcher's works do not yield any positive scientific results. There are many formal manipulations and computations in them; essential difficulties are usually dismissed with a few words without deeper treatment. The content and character diverges significantly from modern research.\\

One should also add:\\
1. The shortcoming, or rather lack of rigor of the definition of iteration with an arbitrary exponent introduced by the candidate met with justified and clearly written criticism by Dr. Stanis\l aw Ruziewicz in {\it Wektor}, Warsaw 1912, no. 5 [On a problem concerning commuting functions].\\
2. Dr. B\"ottcher applies for a second time for veniam legendi in mathematics. The first time the candidate was advised to withdraw his application because of the faults that the Committee at that time found with the candidate's works. These  faults and inadequacies were of the same nature which characterizes the candidate's work also today.\\

The Committee's decision  passed unanimously on June 21, 1918: Not to admit Dr. B\"ottcher to further stages of habilitation. (signature illegible).\\


\begin{center}

Attempt 4, 1919

\includegraphics[height=15cm]{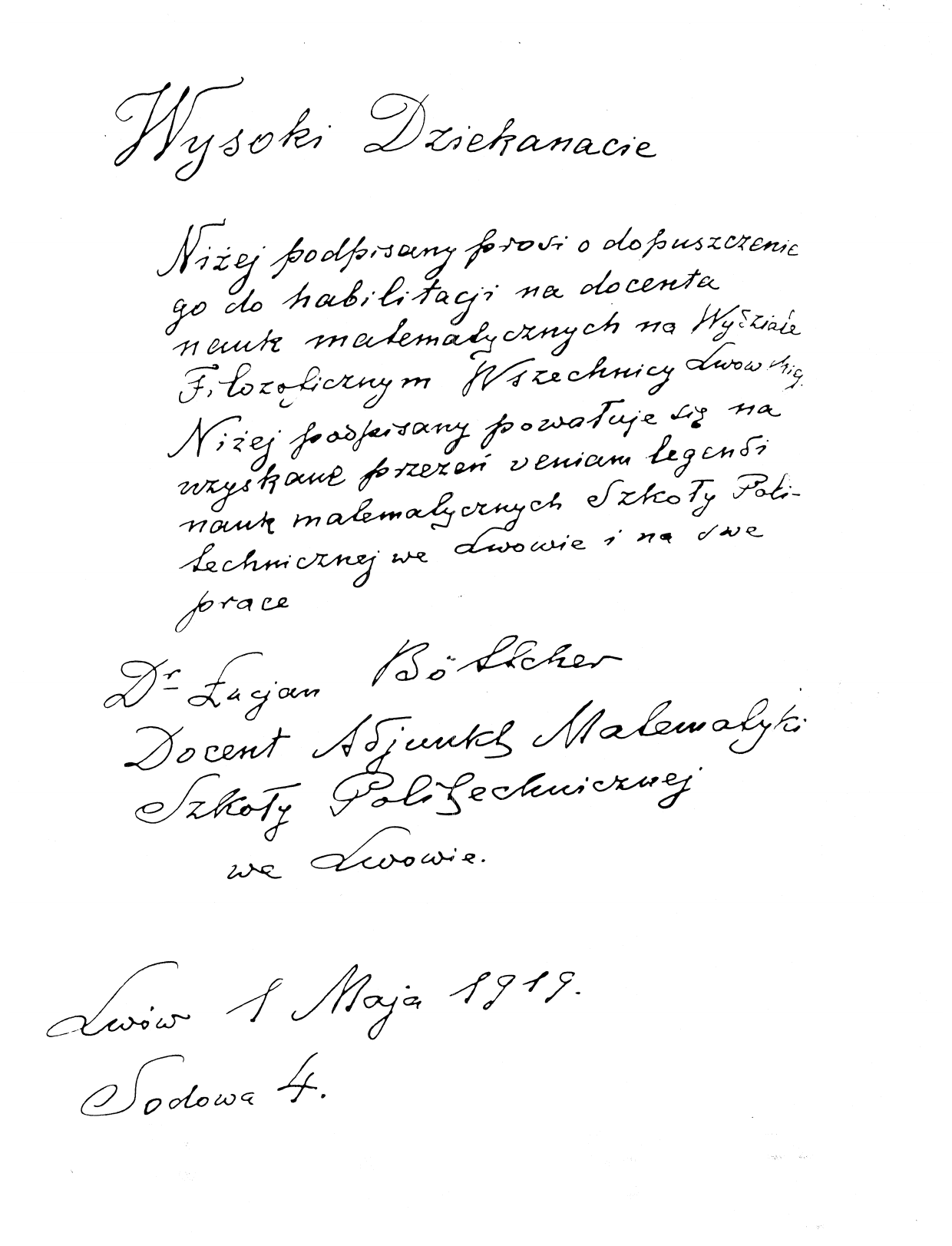}

\end{center}


Honourable Dean's Office\\

The undersigned applies for admission to habilitation for a docent of mathematical sciences in the Faculty of Philosophy of the Lvov University. The undersigned refers to his obtaining veniam legendi in mathematical sciences in the Polytechnic School in Lvov as well as to two papers.\\

Dr. Lucjan B\"ottcher\\
Docent Adiunkt in Mathematics\\
of the Polytechnic School\\
in Lvov\\

Lvov, May 1, 1919, Sodowa 4\\

\subsection{Programmes of some of L. B\"ottcher's lectures at the\\ c.k. Polytechnic School in Lvov (known after 1918 as the Lvov Polytechnic)}

\textbf{METHODS OF COMPUTATION}\\

Graphical way of solving systems of two or three equations with the same number of unknowns. Computing values of technically most significant power, logarithmic, goniometric and cyclometric expressions. Logarithmic and goniometric way of solving equations of second and third degree.\\

\textbf{THEORY OF DIFFERENCE EQUATIONS}\\

Differential equations and difference equations. Solving elementary difference equations. Linear difference equations of order one and higher. Technical applications.\\

\textbf{THEORY OF VECTORS}\\

Development of the notions of a scalar and vector quantity in their geometrical and arithmetical aspects. Principles of the vector calculus and their application to major problems of mechanics, physics and electrotechnics.\\

\textbf{APPLIED MATHEMATICS}\\

Making graphs of most technically significant functions. Reading ready-made graphs. Functional scale and logarithmic sliding rule. Principles of nomography. Graphical way of  solving algebraic and differential equations.\\

B\"ottcher's registry card from the Lvov Polytechnics (including the information about his family)

\begin{center}

\includegraphics[height=15cm]{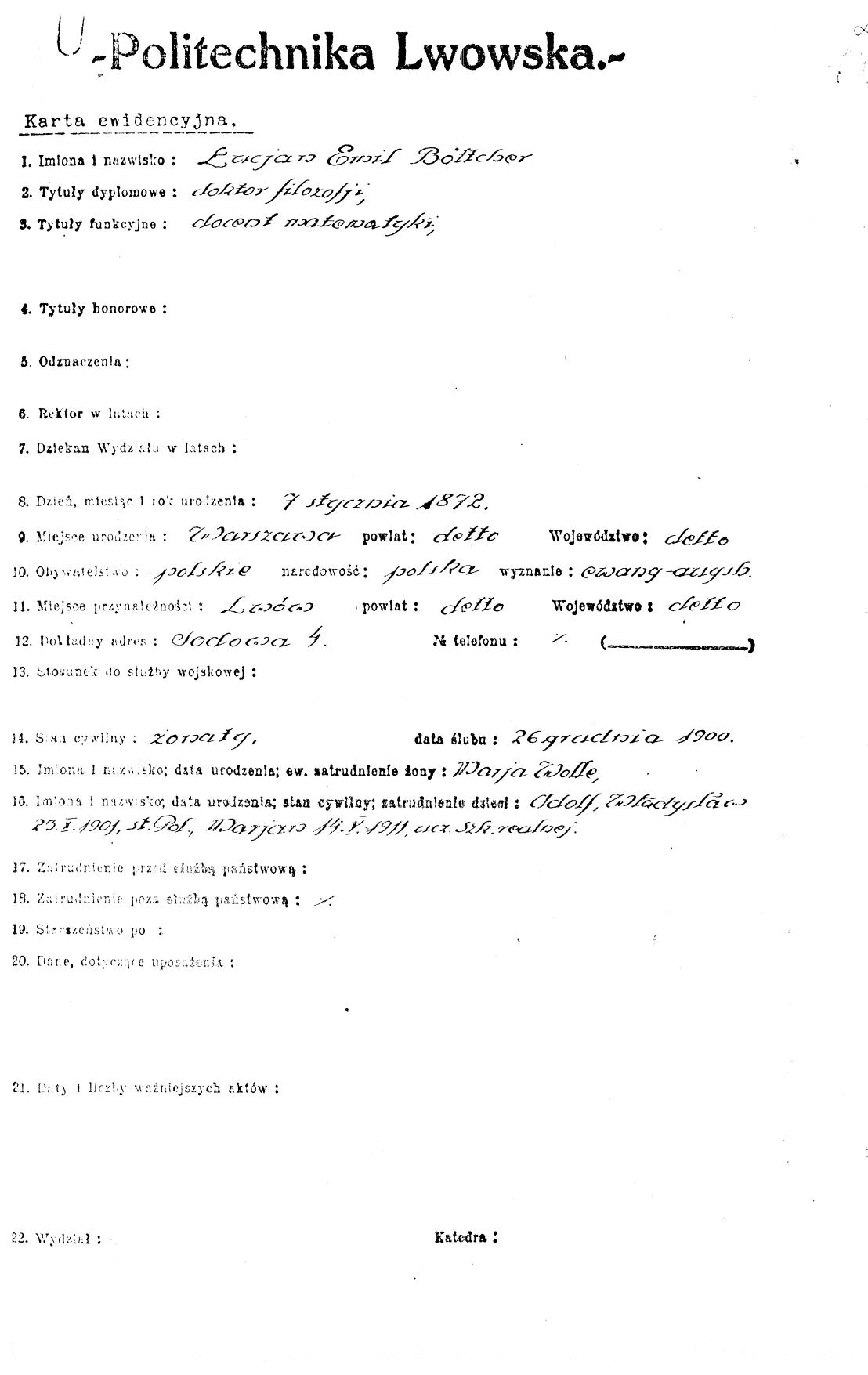}

\end{center}

\newpage

A current (2012) photograph of the house in Sodowa 4 in Lvov where B\"ottcher lived (taken by S. Domoradzki)

\begin{center}

\includegraphics{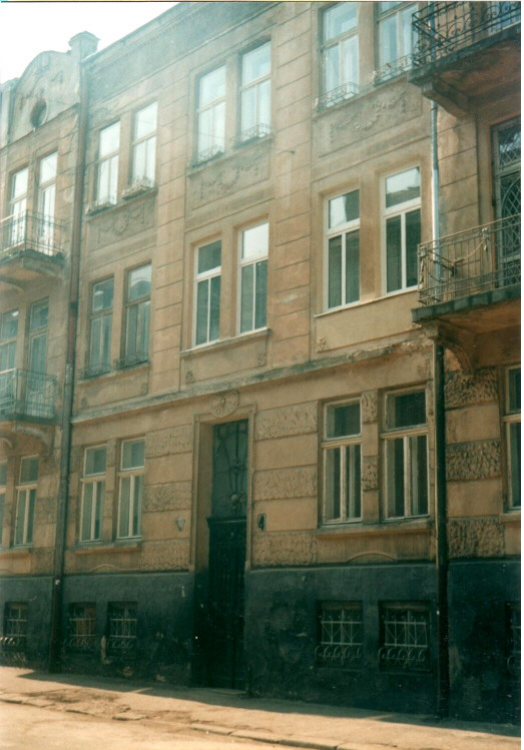}

\end{center}









\end{document}